\documentclass[graybox]{svmult}

\usepackage{type1cm}        
\usepackage{makeidx}         
\usepackage{graphicx}        
\usepackage{multicol}        
\usepackage[bottom]{footmisc}

\usepackage{newtxtext}       %
\usepackage{newtxmath}       

\makeindex             

\newcommand{\R}{\mathbb{R}}
\newcommand{\C}{\mathbb{C}}

\newcommand{\Z}{\mathbb{Z}}

\newcommand{\LL}{\mathrm{L}}

\newcommand{\hyper}[5]
{{}_{#1} F_{#2} \!\left[
  \begin{array}{l}
    #3;\\#4;
  \end{array}#5\right]
}
\newcommand{\HH}{\mathrm{H}}

\newcommand{\PW}{\mathrm{PW}}

\newenvironment{Eqnarray*}
{\arraycolsep=1.4pt
  \begin{eqnarray*}}
  {\end{eqnarray*}
  \hspace*{-4pt}}

\begin{document}
  
   \title*{A Differential Analogue of Favard's Theorem}
  \author{Arieh Iserles and Marcus Webb}
  \institute{Arieh Iserles \at Department of Applied Mathematics and Theoretical Physics, Centre for Mathematical Sciences, University of Cambridge, Wilberforce Rd, Cambridge CB4 1LE, United Kingdom, \email{ai10@cam.ac.uk}
    \and Marcus Webb \at Department of Mathematics, University of Manchester, Alan Turing Building, Manchester M13 9PL, United Kingdom, \email{marcus.webb@manchester.ac.uk}}

  \maketitle
  
  \abstract{Favard's theorem characterizes bases of functions $\{p_n\}_{n\in\Z_+}$ for which $x p_n(x)$ is a linear combination of $p_{n-1}(x)$, $p_n(x)$, and $p_{n+1}(x)$ for all $n \geq 0$ with $p_{0}\equiv1$ (and $p_{-1}\equiv 0$ by convention). In this paper we explore the differential analogue of this theorem, that is, bases of functions $\{\varphi_n\}_{n\in\Z_+}$ for which $\varphi_n'(x)$ is a linear combination of $\varphi_{n-1}(x)$, $\varphi_n(x)$, and $\varphi_{n+1}(x)$ for all $n \geq 0$ with $\varphi_{0}(x)$ given (and $\varphi_{-1}\equiv 0$ by convention). We answer questions about orthogonality and completeness of such functions, provide characterisation results, and also, of course, give plenty of examples and list challenges for further research. Motivation for this work originated in the numerical solution of differential equations, in particular spectral methods which give rise to highly structured matrices and stable-by-design methods for  partial differential equations of evolution. However, we believe this theory to be of interest in its own right, due to the interesting links between orthogonal polynomials, Fourier analysis and Paley--Wiener spaces, and the resulting identities between different families of special functions.}

  \section{Introduction}
  
  Favard's theorem \cite[p.~21]{chihara78iop} states that a sequence of univariate functions $P = \{p_n\}_{n\in\Z_+}$ satisfies
  \begin{eqnarray*}
  p_0(x) &=& 1 \\
  xp_0(x) &=& \beta_0 p_0(x) + \gamma_0p_1(x) \\
  xp_n(x) &=& \alpha_n p_{n-1}(x) + \beta_n p_n(x) + \gamma_np_{n+1}(x), \qquad n \geq 1,
  \end{eqnarray*}
  for some real sequences $\{\alpha_n\}_{n=1}^\infty$, $\{\beta_n\}_{n\in\Z_+}$, $\{\gamma_n\}_{n\in\Z_+}$ if and only if $p_n$ is a polynomial of degree $n$ for all $n \geq 0$ and there exists a finite, signed Borel measure $\mu$ such that
  \begin{eqnarray*}
    \int_{-\infty}^\infty \,\D\mu(x) &=& 1, \\
  \int_{-\infty}^\infty p_m(x) p_n(x) \, \D \mu(x) &=& 0 \text{ for all } n \neq m.
  \end{eqnarray*}
  In other words, $P$ are orthogonal polynomials for the measure $\mu$. Favard's theorem was first announced by Favard in 1935 \cite{favard35sur}, but was discovered independently at about the same time by Shohat and Natanson.
  
  The measure $\mu$ is clearly only unique up to a diagonal scaling of the polynomial sequence $P$, but it is also possible that $\mu$ is non-unique if there exists another measure with the same moments $M_n = \int_{-\infty}^\infty x^n \, \D\mu(x)$ (i.e. the moment problem is indeterminate \cite{aheizer1965classical}). 
  
 The following observations on Favard's theorem can be readily checked \cite{chihara78iop}.
  \begin{itemize}
    \item $\mu$ is a positive measure if and only if $\alpha_n\gamma_{n-1} > 0$ for all $n$.
    \item $\mu$ is a symmetric measure, i.e. $\D\mu(-x) = \D\mu(x)$, if and only if $\beta_n = 0$ for all $n$.
    \item $P$ are orthonormal polynomials with respect to a symmetric measure $\mu$ if and only if $\alpha_n = \gamma_{n-1}$ for $n \geq 1$.
    \end{itemize}

  The topic of this paper is the \emph{differential analogue} of Favard's theorem, and its consequences. Specifically, this is  the characterization of sequences of functions $\Phi = \{\varphi_n\}_{n\in\Z_+}$ which satisfy
 \begin{eqnarray*}
   \varphi_0'(x) &=& \beta_0 \varphi_0(x) + \gamma_0\varphi_1(x) \\
   \varphi_n'(x) &=& \alpha_n \varphi_{n-1}(x) + \beta_n \varphi_n(x) + \gamma_n\varphi_{n+1}(x), \qquad n \geq 1,
 \end{eqnarray*}
 for some real sequences $\{\alpha_n\}_{n=1}^\infty$, $\{\beta_n\}_{n\in\Z_+}$, $\{\gamma_n\}_{n\in\Z_+}$, where $\varphi_0$ is taken as a given smooth function. Despite the fact that differentiation is arguably a purely real-valued affair, we will see that complex coefficients can arise naturally.
  
  The most well known example which satisfies such a differential recurrence is the set of Hermite functions, familiar in mathematical physics:
  \begin{equation}
    \label{HF1}
    \varphi_n(x)=\frac{(-1)^n}{(2^nn!)^{1/2}\pi^{1/4}} \E^{-x^2/2}\mathrm{H}_n(x),\qquad n\in\Z_+,\quad x\in\R,
  \end{equation}
  where $\mathrm{H}_n$ is the $n$th {\em Hermite polynomial.} It is known that
  \begin{equation}
    \label{HF_phi}
    \varphi_n'(x)=-\sqrt{\frac{n}{2}}\varphi_{n-1}(x)+\sqrt{\frac{n+1}{2}}\varphi_{n+1}(x),\qquad n\in\Z_+.
  \end{equation}
  
  Motivation for studying such bases is numerical solution of differential equations. If a solution is approximated using a \emph{spectral method} in the basis $\Phi$, numerical solutions take the form $u_N(x) = \sum_{n=0}^N a_n \varphi_n(x)$, and the derivative (after projecting back onto the span of the first $N+1$ basis elements) is given by $u_N'(x) = \sum_{n=0}^{N} \tilde{a}_n \varphi_n(x)$, where $\tilde{\mathbf{a}} = D \mathbf{a}$, and
  \begin{equation}
  \label{diffmat}
  D = \left( \begin{array}{ccccc}
               \beta_0  & \alpha_1 &              &              &          \\
               \gamma_0 & \beta_1  & \alpha_2     &              &          \\
                        & \gamma_1 & \beta_2      & \ddots       &          \\
                        &          & \ddots       & \ddots       & \alpha_N \\
                        &          &              & \gamma_{N-1} & \beta_N
               \end{array} \right)
  \end{equation}
 The matrix $D$ is called the \emph{differentiation} matrix for the basis $\Phi$. Sparse and structured matrices such as this can be used to an advantage for the numerical linear algebra involved in solving a differential equation numerically.
 
 Beyond the tridiagonal nature of $D$ which is guaranteed by such an analogue of Favard's theorem, we have the following wish list:
 \begin{itemize}
   \item We wish for $D$ to be skew-symmetric (i.e. $\gamma_n = -\alpha_{n+1}$). This emulates the skew-self-adjointness of the differentiation operator with respect to the standard inner product and zero Dirichlet (or Cauchy) boundary conditions. For numerical solution of time-dependent PDEs this ensures stability in the $\ell_2$ norm on the expansion coefficients.  Moreover, solutions of numerous dispersive PDEs, e.g.\ the linear Schr\"odinger equation, preserve the $\mathrm{L}_2$ norm, and skew-self-adjointness of $D$ ensures that their discretisation does so in the $\ell_2$ norm as well.
   \item We wish for $\Phi$ to be complete in the sense that all functions of interest can be approximated by finite linear combinations in the basis.
   \item We wish for $\Phi$ to be an orthonormal basis, since then the discrete $\ell_2$ norm on the expansion coefficients is precisely equal to the continuous $\mathrm{L}_2$ norm of the numerical solution, and so the stability properties transfer over to the continuous norm.
   \item We wish for $\varphi_n(x)$ to be easily computable, perhaps with an explicit analytical formula, and for expansions like $u_N(x)$ to be easy to work with computationally.
 \end{itemize}

Most of the results in this paper have been published elsewhere by the present authors \cite{iserles19oss,iserles20for,iserles20fast}. Our purposes here are to provide these results in a self-contained fashion, to present the material in a way more suited for experts in classical analysis, and to offer some new examples and insight that generalize the theory.

  \section{The main theory}
 
At a basic level, the differential analogue of Favard's theorem is solved by considering functions of the form,
 \begin{equation*}
   \varphi_n(x) = p_n\left( \frac{\D }{\D x}\right) \varphi_0(x),
   \end{equation*}
 where $\{p_n\}_{n\in\Z_+}$ are the orthogonal polynomials with recurrence coefficients $\{\alpha_n\}_{n=1}^\infty$, $\{\beta_n\}_{n\in\Z_+}$, $\{\gamma_n\}_{n\in\Z_+}$. However, this formulation does not help us design the differentiation matrix to be skew-symmetric unless we allow non-positive measures $\mu$, nor does it elucidate the function space spanned by $\{\varphi_n\}_{n\in\Z_+}$ or when the basis will be orthogonal. Cue the following three theorems and their corollaries, which approach the problem via Fourier analysis \cite{iserles19oss,iserles20for}.
  
  \begin{theorem}\label{thm:nonorthogonalPhi} 
    A set of functions  $\Phi = \{\varphi_n\}_{n\in\Z_+}\subset\mathrm{L}_2(\R)$ satisfies
    \begin{eqnarray*}
      \varphi_0' &=& \I c_0 \varphi_0 + b_0 \varphi_1 \\
    \varphi_{n}' &=& -\overline{b_{n-1}} \varphi_{n-1} + \I c_n \varphi_n + b_n \varphi_n 
    \end{eqnarray*}
    for some sequences $b_n \in \C \setminus \{0\}$, $c_n \in \R$, if and only if
    \begin{equation}\label{explicit_phi}
    \varphi_n(x) = \frac{\E^{\I \theta_n}}{\sqrt{2\pi}} \int_{-\infty}^\infty \E^{\I x \xi} p_n(\xi) g(\xi) \, \D \xi,
    \end{equation}
  where
    \begin{itemize}
      \item $P = \{p_n\}_{n\in\Z_+}$ is an orthonormal polynomial system with respect to $\D\mu$, a probability measure on the real line with all moments finite and with infinitely many points of increase.
      \item $g \in \mathrm{L}^2(\R)$ is such that $|\xi^n g(\xi)| \to 0$ as $|\xi| \to \infty$ for all $n \geq 0$
      \item $\{\theta_n\}_{n\in\Z_+} \subset [0,2\pi)$ 
    \end{itemize}
%
  \end{theorem}

\begin{corollary}\label{cor:real}
  A set of real-valued functions $\Phi = \{\varphi_n\}_{n\in\Z_+}\subset\mathrm{L}_2(\R)$ satisfies
  \begin{eqnarray*}
    \varphi_0' &=&  b_0 \varphi_1 \\
    \varphi_{n}' &=& -b_{n-1} \varphi_{n-1} + b_n \varphi_n 
  \end{eqnarray*}
  for some sequence $b_n > 0$, if and only if $\Phi$ satisfies the requirements of Theorem \ref{thm:nonorthogonalPhi} but with the extra constraints,
  \begin{enumerate}
    \item The measure $\mu$ is symmetric, i.e. $\D\mu(-\xi) = \D\mu(\xi)$
    \item The function $g$ has even real part and odd imaginary part
    \item $\theta_n = n\pi/2\bmod {2\pi}$, i.e.\ $\E^{\I\theta_n} = \I^n$.
    \end{enumerate}
  \end{corollary}
  
  \begin{theorem}[Orthogonal systems]\label{thm:orthogonalPhi}
    Let $\Phi = \{\varphi_n\}_{n\in\Z_+}$ satisfy the requirements of Theorem \ref{thm:nonorthogonalPhi}. Then $\Phi$ is orthogonal in $\mathrm{L}_2(\R)$ if and only if $P$ is orthogonal with respect to the measure $|g(\xi)|^2\mathrm{d}\xi$. Furthermore, whenever $\Phi$ is orthogonal, the functions $\varphi_n / \|g\|_{2}$ are orthonormal.
    \begin{proof}
      It follows from Parseval's identity for the inner product of Fourier transforms that
      \begin{equation*}
        \int_{-\infty}^\infty \overline{\varphi_n(x)} \varphi_m(x) \, \D x = \E^{\I(\theta_m-\theta_n)}\int_{-\infty}^\infty p_m(\xi)p_n(\xi)|g(\xi)|^2 \,\D\xi, \qquad n,m\in\mathbb{Z}_+.
        \end{equation*}
      \end{proof}
  \end{theorem}
  
  \begin{theorem}[Orthogonal bases for a Paley--Wiener space]\label{thm:PW}
    Let $\Phi = \{\varphi_n\}_{n\in\Z_+}$ satisfy the requirements of Theorem \ref{thm:orthogonalPhi} with a measure $\mathrm{d}\mu$ such that polynomials are dense in $\mathrm{L}_2(\R;\mathrm{d}\mu)$. Then $\Phi$ forms an orthogonal basis for the Paley--Wiener space $\PW_\Omega(\R)$\footnote{$\PW_\Omega(\R)$ for $\Omega \subseteq \R$ is the space of all functions $f \in \mathrm{L}_2(\R)$ such that the Fourier transform of $f$ is supported on $\Omega$.}, where $\Omega$ is the support of $\mathrm{d}\mu$.
  \end{theorem}
  
  The key consequence of Theorem \ref{thm:PW} is that for a basis $\Phi$ satisfying the requirements of Theorem \ref{thm:orthogonalPhi} to be complete in $\mathrm{L}_2(\R)$, it is necessary that the polynomial basis $P$ is orthogonal with respect to a measure which is supported on the whole real line.

  \section{Examples}\label{sec:examples}
  
  For an absolutely continuous measure $\D\mu(\xi) = w(\xi)\D\xi$ with orthonormal polynomials $P = \{p_n\}_{n\in\Z_+}$, we consider the following to be its canonical associated orthonormal basis functions:
  \begin{equation}\label{eqn:canonical}
    \varphi_n(x) = \frac{\I^n}{\sqrt{2\pi}} \int_{-\infty}^\infty \E^{\I x \xi} p_n(\xi) |w(\xi)|^{\tfrac12} \,\D\xi.
  \end{equation}
  With this particular representation, the coefficients $\{b_n\}_{n\in\Z_+}$ have the convenient property that they are all positive real numbers. The coefficients in the differentiation matrix for $\Phi$ also correspond nicely with the three-term recurrence coefficients of the orthonormal polynomials:
  \begin{eqnarray*}
    \xi p_n(\xi) &=& b_{n-1}p_{n-1}(\xi) + c_n p_n(\xi) + b_n p_{n+1}(\xi) \\
    &\iff& \\
    \varphi_n'(\xi) &=& -b_{n-1}\varphi_{n-1}(\xi) + \I c_n \varphi_n(\xi) + b_n \varphi_{n+1}(\xi)
  \end{eqnarray*}
  
  The terminology we use for these bases is as follows. If the polynomial basis $P$ has a name, such as the \emph{Chebyshev polynomials}, then we call the basis $\Phi$ the \emph{transformed Chebyshev functions}.
  
  \subsection{Jacobi}
  
\begin{figure}[h!]
  \begin{center}
    \includegraphics[width=220pt]{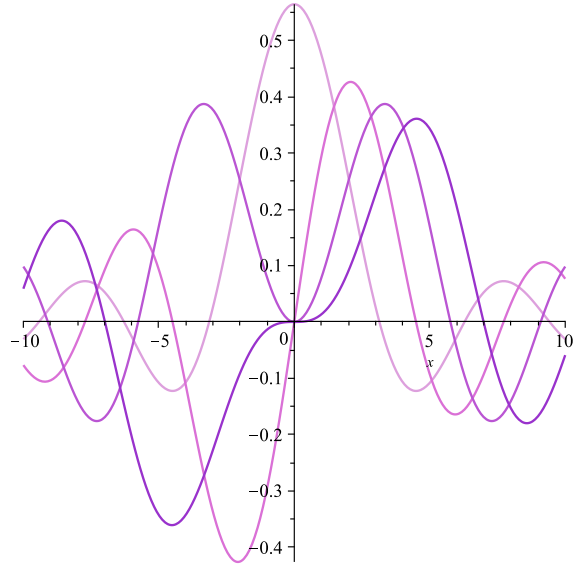}
    \caption{The transformed Legendre functions $\varphi_n$ for $n=0,1,2,3$: darker shade corresponds to higher $n$.}
    \label{Fig:Legendre}
  \end{center}
\end{figure}

Transformed Jacobi functions are related to Bessel functions. The measure of orthogonality is
\begin{equation}\label{eqn:Jacobimeasure}
  \D\mu(\xi) = \chi_{(-1,1)}(\xi) (1-\xi)^\alpha(1+\xi)^\beta \,\D \xi,
\end{equation}
for given $\alpha,\beta > -1$. It is possible to give an analytic formula for the first few transformed basis functions in terms of hypergeometric functions, for example,
\begin{eqnarray*}
  \varphi_0(x) &=& \frac{1}{\sqrt{2\pi}}\int_{-1}^1 (1-\xi)^{\alpha/2}(1+\xi)^{\beta/2} \frac{\E^{\I x\xi}}{\sqrt{2^{\alpha/2+\beta/2+1}B(\alpha/2+1,\beta/2+1)}} \D\xi \\
  &=&  \frac{2^{\alpha/4+\beta/4}}{\sqrt{\pi}}\sqrt{B(\alpha/2+1,\beta/2+1)} \E^{\I x} \hyper{1}{1}{1+\alpha/2}{2+(\alpha+\beta)/2}{-2\I x}.
  \end{eqnarray*}
 However, we do not believe that any traction can be gained from working out these functions for general $\alpha$ and $\beta$. A special case of the Jacobi polynomials is the family of ultraspherical polynomials, $\beta=\alpha$, whereby 
\begin{displaymath}
  b_n=\sqrt{\frac{(n+1)(n+2\alpha+1)}{(2n+2\alpha+1)(2n+2\alpha+3)}},\qquad n\in\mathbb{Z}_+,
\end{displaymath}
and $c_n = 0$ for all $n \in \mathbb{Z}_+$. The first two transformed ultraspherical functions are
\begin{eqnarray*}
  \varphi_0(x)&=&  C_\alpha \left( \frac{2}{x}\right)^{(1+\alpha)/2}\mathrm{J}_{(1+\alpha)/2}(x) \\
  \varphi_1(x) &=& C_\alpha\sqrt{3+2\alpha} \left( \frac{2}{x}\right)^{(1+\alpha)/2}\mathrm{J}_{(3+\alpha)/2}(x),
\end{eqnarray*}
where $\mathrm{J}_\nu(x)$ is the Bessel function of order $\nu > -1$, and 
\begin{displaymath}
  C_\alpha^2 = \frac{\Gamma(\alpha + 3/2)\Gamma(\alpha/2 + 1)}{2^{\alpha+1}\Gamma((1+\alpha)/2)}.
  \end{displaymath}
For higher indices these functions become more complicated expressions involving polynomials in $1/x$ and Bessel functions. However, in the special case $\alpha=0$ of Legendre polynomials these reduce to a very neat form \cite{iserles19oss},
\begin{equation}
  \label{eqn:Legendre}
  \varphi_n(x)=\sqrt{\frac{n+\frac12}{x}} \mathrm{J}_{n+1/2}(x),\qquad n\in\mathbb{Z}_+.
\end{equation}
Fig.~\ref{Fig:Legendre} displays the first four transformed Legendre functions. 

Of course, all functions $\Phi$ `seeded' by Jacobi polynomials are necessarily band limited and live in PW${}_{(-1,1)}(\mathbb{R})$. This could be useful in some applications, such as signal processing, but as stated in the introduction, we would like to generate bases which are complete in $\mathrm{L}_2(\R)$. By Theorem \ref{thm:PW}, we will need to take polynomials which are orthonormal on the whole real line order to achieve this.

\subsection{Hermite}
  
Under the transform in equation \eqref{eqn:canonical}, the Hermite polynomials map directly to Hermite functions, so these polynomials hold a special position as a sort of fixed point in this theory. We already mentioned Hermite functions, generated by $\D\mu(\xi)=\E^{-\xi^2}\D\xi$, in \eqref{HF_phi}. They are widely used in computation on the real line and are a univariate case of Hagedorn wave packets \cite{lasser2020cqd}. It has been proved in \cite{iserles19oss} that they are unique among all systems $\Phi$ consistent with Theorems~2--3 with the representation $\varphi_n(x)=h(x) q_n(x)$, where $h\in\mathrm{L}_2(\mathbb{R})$ and each $q_n$ is a polynomial of degree $n$. 

Hermite functions obey the Cram\'er inequality $|\varphi_n(x)|\leq\pi^{-1/4}$, $x\in\mathbb{R}$, $n\in\mathbb{Z}_+$, and several helpful identities, while their generating function is inherited in a straightforward manner from Hermite polynomials. 

\begin{figure}[h!]
  \begin{center}
    \includegraphics[width=220pt]{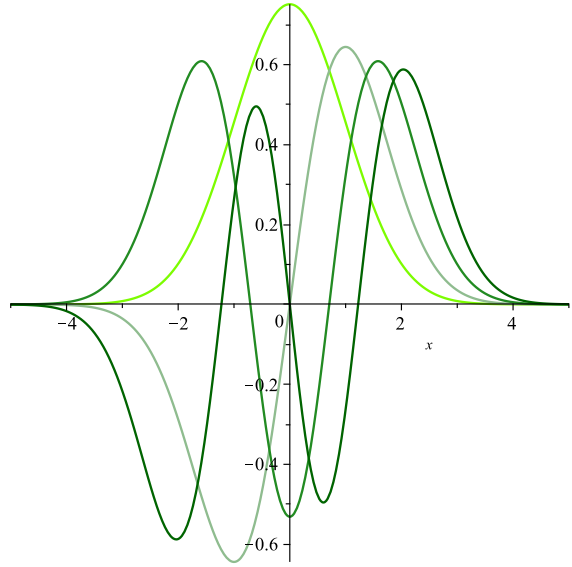}
    \caption{Hermite functions $\varphi_n$ for $n=0,1,2,3$: darker shade corresponds to higher $n$.}
    \label{Fig:Hermite}
  \end{center}
\end{figure}

\begin{figure}[tb]
  \begin{center}
    \includegraphics[width=220pt]{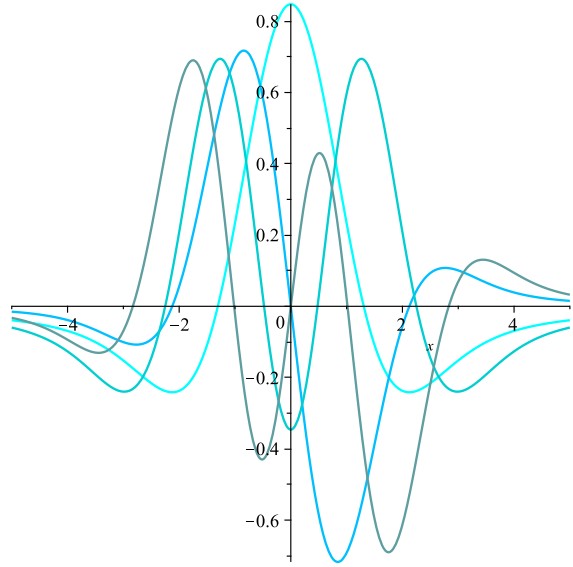}
    \caption{Transformed generalized Hermite functions ($\eta=1$) $\varphi_n$ for $n=0,1,2,3$: darker shade corresponds to higher $n$.}
    \label{Fig:GenHermite}
  \end{center}
\end{figure}
  
  \subsection{Generalized Hermite}
  
  Let $\eta>-1/2$. The generalized Hermite polynomials $\HH_n^{(\eta)}$ are orthogonal with respect to
  \begin{displaymath}
    \D\mu(\xi)=|\xi|^{2\eta}\E^{-\xi^2}\D\xi,\qquad \xi\in\mathrm{R}
  \end{displaymath}
  and obey the three-term recurrence relation
  \begin{displaymath}
    \HH_{n+1}^{(\eta)}(\xi)=2\xi\HH^{(\eta)}_n(\xi)-2(n+\theta_n)\HH_{n-1}^{(\eta)}(\xi),
  \end{displaymath}
  where $\theta_n=0$ for an even $n$ and $\theta_n=2\eta$ otherwise \cite[p,~156--157]{chihara78iop}. The explicit form of $\Phi$ has been derived in \cite{iserles19oss} -- the algebra is fairly laborious. Simplifying somewhat the formul\ae{} therein with the first Kummer formula \cite[p.~125]{rainville60sf}, we have
  \begin{eqnarray*}
    \varphi_{2n}(x)&=&c_0\frac{(-1)^n(\frac12+\frac{\eta}{2})_n}{2^n}\!\sum_{\ell=0}^n {n\choose\ell} \frac{(-1)^\ell(-n\!-\!\eta\!+\!\frac12)_\ell}{2^\ell(-n\!+\!\frac12\!-\!\frac{\eta}{2})_\ell} \hyper{1}{1}{n\!-\!\ell\!+\!\frac{\eta+1}{2}}{\frac12}{-\frac{x^2}{2}},\\
    \varphi_{2n+1}(x)&=&c_0\frac{(-1)^n(\frac32+\frac{\eta}{2})_n}{2^n}x \!\sum_{\ell=0}^n{n\choose\ell} \frac{(-1)^\ell(-n\!-\!\eta\!-\!\frac12)_\ell}{2^\ell(-n\!-\!\frac12\eta\!-\!\frac12)_\ell} \hyper{1}{1}{n\!-\!\ell\!+\!\frac{\eta+3}{2}}{\frac32}{-\frac{x^2}{2}},
  \end{eqnarray*}
  where the normalising constant $c_0$ can be obtained from $c_0^2\int_{-\infty}^\infty \varphi_0^2(x)\D x=1$, $c_0>0$.
  
Figs \ref{Fig:Hermite}--\ref{Fig:GenHermite} display Hermite functions and transformed generalized Hermite system for $\eta=1$.  Hermite functions are, needless to say, the familiar Hermite polynomials scaled by $\E^{-x^2}$ and they demonstrate rapid decay, while their zeros are all real and interlace. Less is known about the case $\eta=1$ except that the functions evidently decay more sedately. 
 
  \subsection{Laguerre}
  
  The transformed Laguerre functions are related to the Fourier basis. The Laguerre weight is
  \begin{equation*}
    \D\mu\xi = \chi_{[0,\infty)}(\xi) \E^{-\xi} \,\D\xi.
    \end{equation*}
  
  As described in \cite{iserles20for},  transformed Laguerre functions have the particularly elegant form,
  \begin{displaymath}
    \varphi_n(x)=\sqrt{\frac{2}{\pi}} \I^n \frac{(1+2\I x)^n}{(1-2\I x)^{n+1}},\qquad n\in\Z_+.
  \end{displaymath}
  These functions do not form a complete orthonormal basis for $\mathrm{L}_(\R)$; by Theorem \ref{thm:PW} they are dense in $\mathrm{PW}_{[0,\infty)}(\R)$. If we want to obtain a basis which is dense in $\mathrm{L}_2(\R)$, all we need to do is take a direct sum of these functions with the functions associated to the Laguerre measure on $(-\infty,0]$, namely $\D\mu(\xi) = \chi_{(-\infty,0]}(\xi) \E^{\xi} \,\D\xi$. This corresponds to simply taking the above formula and indexing it by $n \in \Z$, resulting in what are known in harmonic analysis as the \emph{Malmquist--Takenaka functions},
  \begin{equation}
    \label{MalmTak}
    \varphi_n(x)=\sqrt{\frac{2}{\pi}} \I^n \frac{(1+2\I x)^n}{(1-2\I x)^{n+1}},\qquad n\in\Z.
  \end{equation}
  These functions have a wealth of beautiful properties and have been discovered and rediscovered over nearly a century since their initial discovery by Malmquist \cite{malmquist26stc} and Takenaka \cite{takenaka25oof}, both in 1926.
  
  The most notable property of the Malmquist--Takenaka basis is its relation to the Fourier basis. If we make the change of variables $\theta = 2\arctan(2 x)$ and $x = \tfrac12 \tan\left(\tfrac12\theta\right)$, then
  \begin{equation}\label{eqn:MTFourier}
    \varphi_n(x) = \sqrt{\frac{2}{\pi}}\I^n\E^{\I (n+\tfrac12) \theta} \cos \frac{\theta}{2}.
    \end{equation}
  So we see that the Malmquist--Takenaka basis is the Fourier basis in disguise, which leads to a fast FFT-based algorithm to compute the expansion coefficients of a given $f \in \mathrm{L}_2(\R)$ (see Section \ref{sec:computational}).
  
  \begin{figure}[tb]
\begin{center}
  \includegraphics[width=100pt]{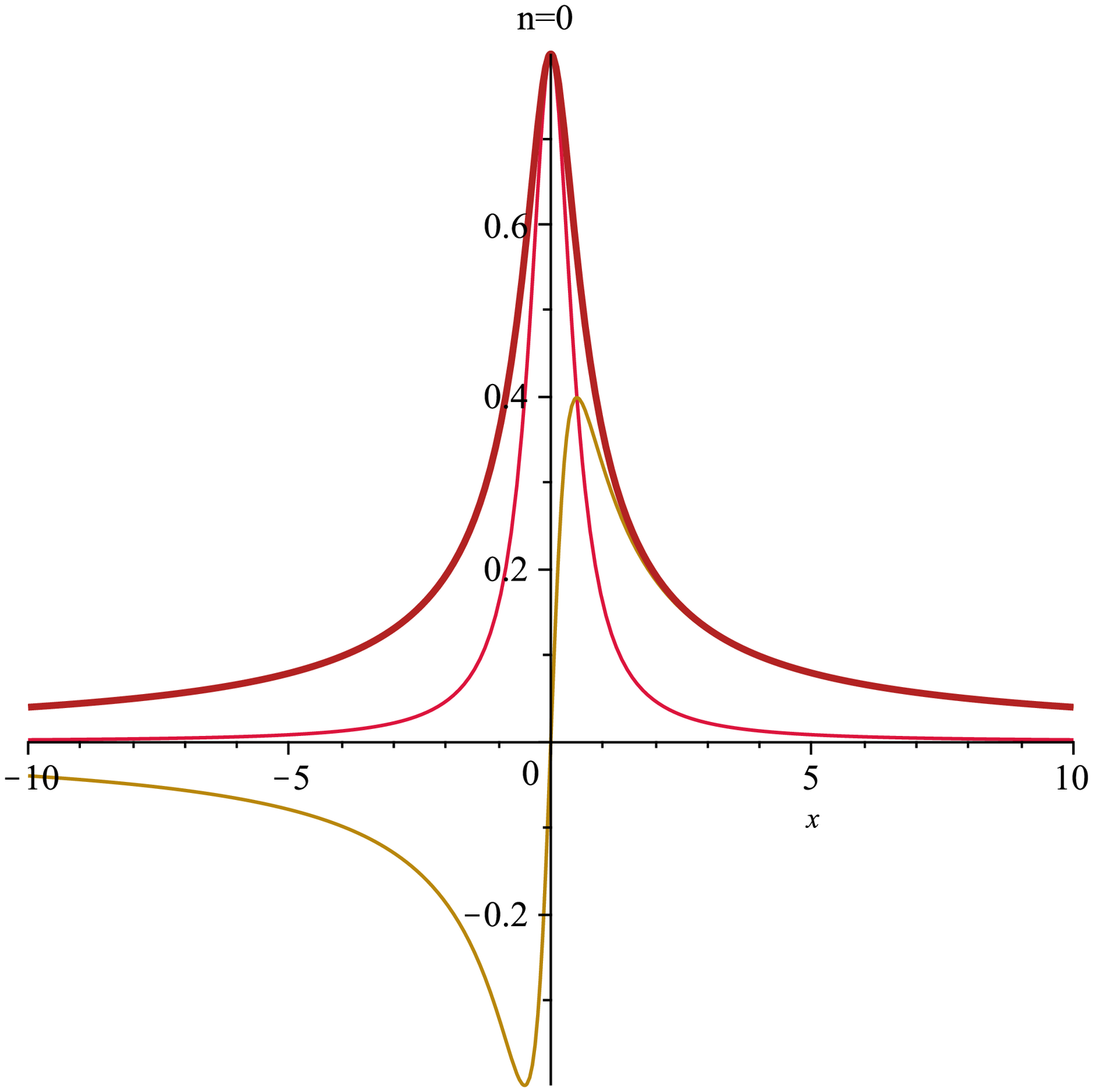}\hspace*{15pt}\includegraphics[width=100pt]{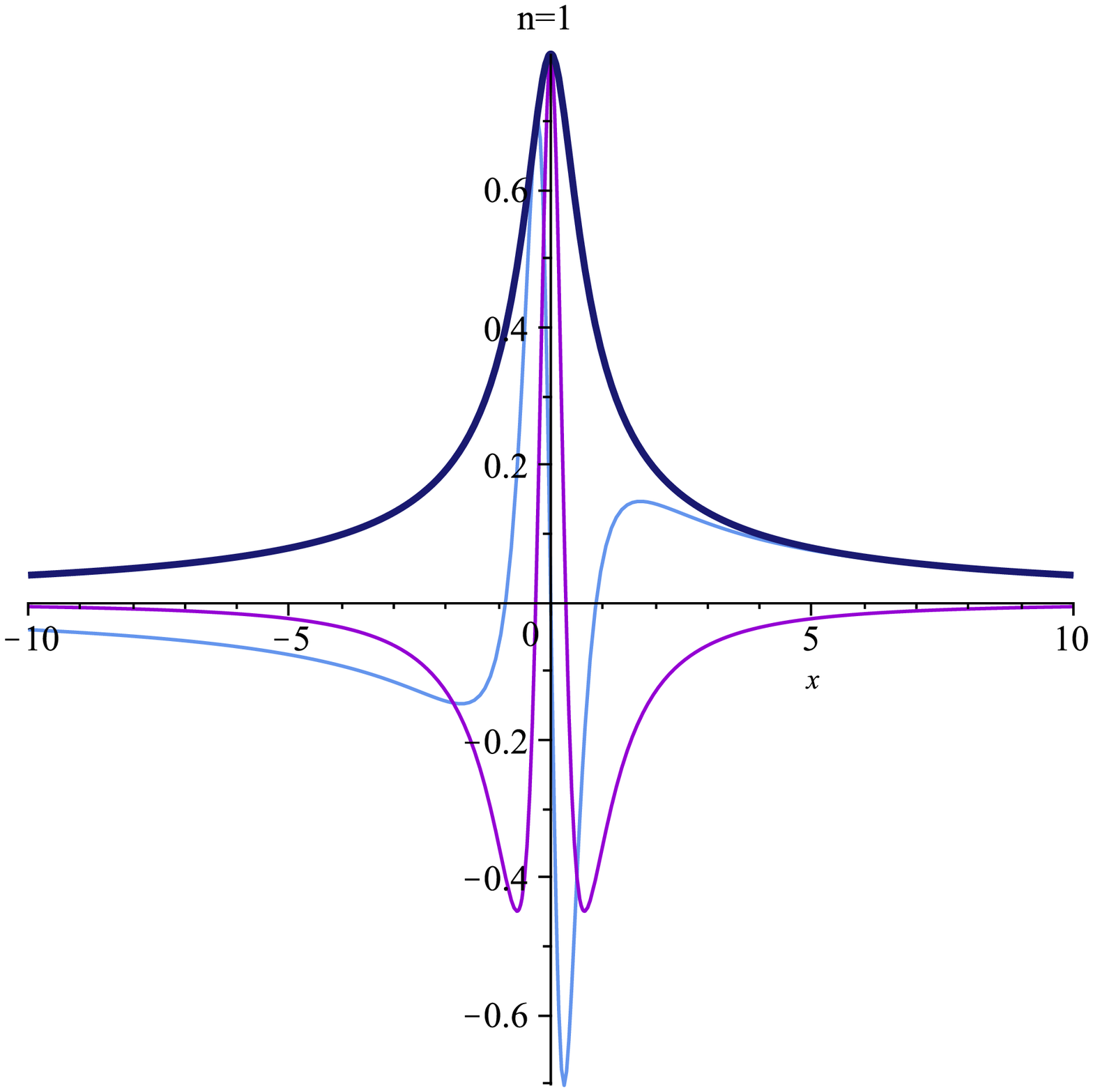}\hspace*{15pt}\includegraphics[width=100pt]{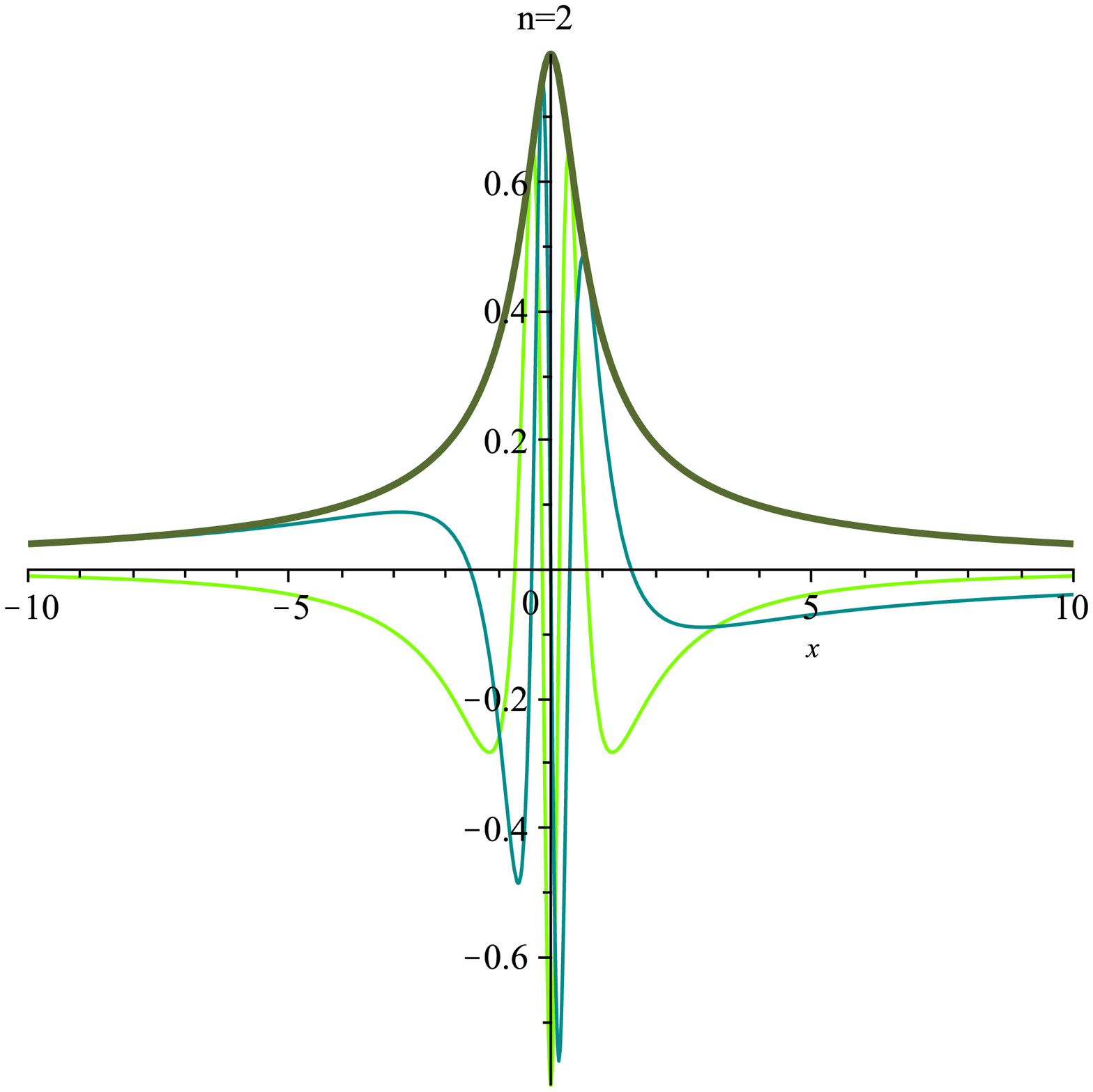}
  
  \includegraphics[width=100pt]{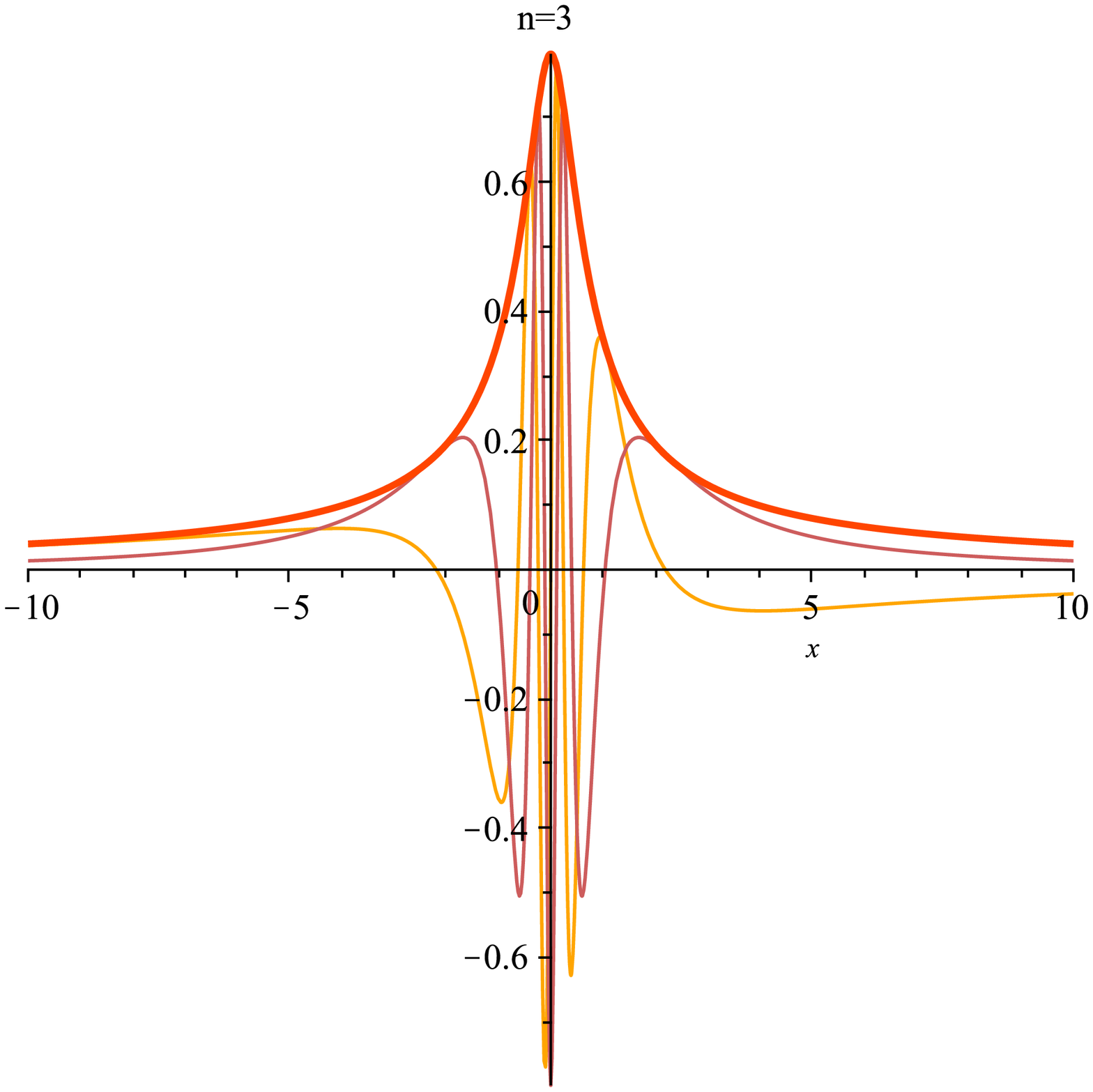}\hspace*{15pt}\includegraphics[width=100pt]{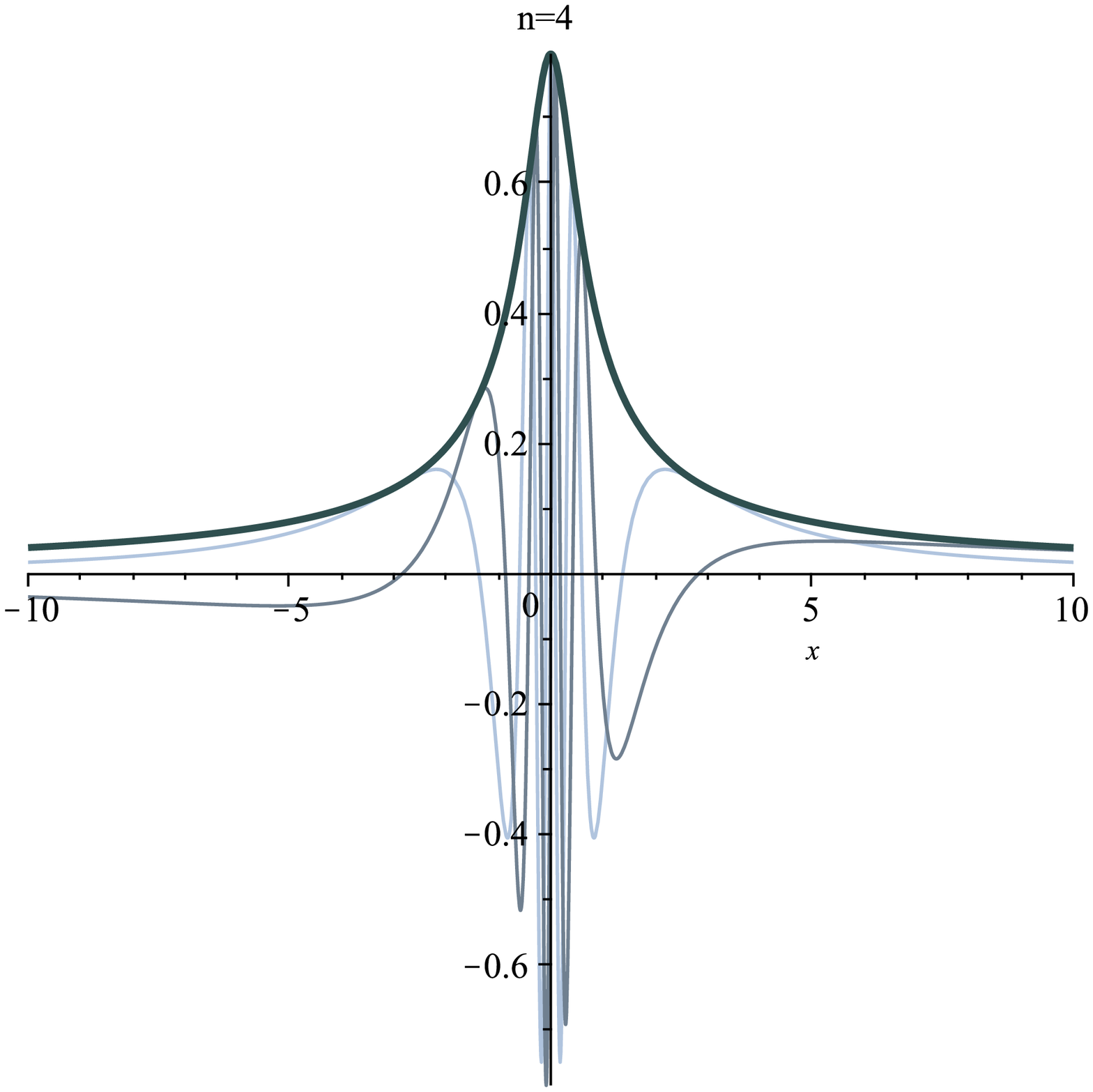}
\caption{The Malmquist--Takenaka system for $n=0,1,2,34$; Real part in lighter shade than the imaginary part, both within the envelope $(1+4x^2)^{-1/2}$, plotted by a thicker line.}
\label{Fig:MT}
\end{center}
\end{figure}

  \subsection{Generalized Laguerre}
  
  The transformed generalized Laguerre functions are related to Szeg\H{o}--Askey polynomials on the unit circle. The {\em generalized Laguerre measure\/} is given by,
  \begin{displaymath}
    \D\mu(\xi)=\chi_{(0,\infty)}(\xi) \xi^\alpha \E^{-\xi}\D \xi,
    \end{displaymath} 
  where $\alpha>-1$. The case $\alpha = 0$ gives the standard Laguerre polynomials. The corresponding orthogonal polynomials are the {\em generalized Laguerre polynomials\/}
  \begin{displaymath}
    \LL_n^{(\alpha)}(\xi)=\frac{(1+\alpha)_n}{n!} \hyper{1}{1}{-n}{1+\alpha}{\xi}=\frac{(1+\alpha)_n}{n!} \sum_{\ell=0}^n (-1)^\ell {n\choose \ell} \frac{\xi^\ell}{(1+\alpha)_\ell},
  \end{displaymath}
  where $(z)_m=z(z+1)\cdots(z+m-1)$ is the {\em Pochhammer symbol\/} and ${}_1F_1$ is a {\em confluent hypergeometric function\/} \cite[p.~200]{rainville60sf}. The Laguerre polynomials obey the recurrence relation
  \begin{displaymath}
    (n+1)\LL_{n+1}^{(\alpha)}(\xi)=(2n+1+\alpha-\xi)\LL_n^{(\alpha)}(\xi)-(n+\alpha)\LL_{n-1}^{(\alpha)}(\xi).
  \end{displaymath}
\begin{displaymath}
  p_n(\xi)=(-1)^n \sqrt{\frac{n!}{\Gamma(n+1+\alpha)}} \LL_n^{(\alpha)}(\xi),\qquad n\in\Z_+. 
\end{displaymath}
We deduce after simple algebra that the normalized polynomials have three-term recurrence coefficients given by,
\begin{displaymath}
  b_n=\sqrt{(n+1)(n+1+\alpha)},\qquad c_n=2n+1+\alpha.
\end{displaymath}
After lengthy algebra, it was shown in \cite{iserles20for} that the transformed generalized Laguerre functions can be expressed by
  \begin{equation}
    \label{FL}
    \varphi_n(x) = (-\I)^n\sqrt{\frac{2}{\pi}} \left(\frac{1}{1-2\I x}\right)^{\!1+\frac{\alpha}{2}} \Pi^{(\alpha)}_n\left(\frac{1+2\I x}{1-2\I x}\right), 
  \end{equation}
  where $\Pi^{(\alpha)}_n$ is a polynomial of degree $n$. Using the substitution $x = \frac12 \tan \frac{\theta}{2}$ for $\theta \in (-\pi,\pi)$, which implies $(1+2\I x)/(1-2\I x) = \E^{\I \theta}$, the orthonormality of the basis $\Phi$ can be seen to imply that $\{\Pi^{(\alpha)}_n\}_{n\in\Z_+}$ are in fact orthogonal polynomials on the unit circle (OPUC) with respect to the weight
  \begin{displaymath}
    W(\theta) = \cos^{\alpha} \frac{\theta}{2}.
  \end{displaymath}
  To be clear, this means that for all $n,m\in\Z_+$,
  \begin{displaymath}
    \frac{1}{2\pi}\int_{-\pi}^\pi \overline{\Pi^{(\alpha)}_n(\E^{\I \theta})}\Pi^{(\alpha)}_m(\E^{\I \theta}) \,\cos^{\alpha} \frac{\theta}{2}\,\mathrm{d} \theta = \delta_{n,m}.
  \end{displaymath}
  These polynomials are related to the \emph{Szeg\H{o}--Askey polynomials} \cite[18.33.13]{DLMF}, $\{\phi^{(\lambda)}_n\}_{n\in\Z_+}$, which satisfy
  \begin{displaymath}
    \frac{1}{2\pi}\int_{-\pi}^\pi \overline{\phi^{(\lambda)}_n(\E^{\I \theta})}\phi^{(\lambda)}_m(\E^{\I \theta}) \,(1-\cos \theta)^\lambda\,\mathrm{d} \theta = \delta_{n,m},
  \end{displaymath}
  by the relation $\Pi^{(\alpha)}_n(z) \propto \phi^{(\alpha/2)}_n(-z)$. In turn, these polynomials are related to the Jacobi polynomials $\mathrm{P}_n^{(\tfrac{\alpha-1}{2},-\tfrac12)}$ and $\mathrm{P}_n^{(\tfrac{\alpha+1}{2},\tfrac12)}$ via the Delsarte--Genin relationship \cite{szego1939orthogonal}.

  \subsection{Continuous Hahn}
    
  Continuous Hahn polynomials \cite{koekoek2010hypergeometric} are orthogonal with respect to the measure
  \begin{displaymath}
     \D\mu_{a,b}(\xi) =  \frac{1}{2\pi}|\Gamma(a+\I \xi)\Gamma(b-\I \xi)|^2\D\xi, \qquad \xi\in (-\infty,\infty),
    \end{displaymath}
  where $a$ and $b$ are complex parameters with positive real part. The transformed continuous Hahn functions are related to Jacobi polynomials in the following way. Take the following square root of the measure,
  \begin{displaymath}
    g_{a,b}(\xi) = \frac{1}{\sqrt{2\pi}}\Gamma(a+\I\xi)\Gamma(b-\I\xi).
  \end{displaymath} 
  Note that in general $g_{a,b}$ is a complex valued function, hence it deviates from what we declared as `canonical' at the beginning of Section \ref{sec:examples}. Note further that when $a$ and $b$ are real, it has an even real part and odd imaginary part, so the resulting transformed functions can be made to be real-valued in this case by Corollary \ref{cor:real}.
  
  In \cite{iserles20fast}, the following remarkable identity was shown (in fact an analogous identity was shown for the non-standard continuous Hahn measure $\D\mu_{a,b}(\xi/2)$). Let $p_n$ be the normalized continuous Hahn polynomials (with parameters $a$ and $b$), then
  \begin{eqnarray}\label{contHahn_varphi}
    &&\frac{\I^n}{2\pi} \int_{-\infty}^\infty p_n(\xi) \Gamma(a+\I\xi)\Gamma(b-\I\xi) \, \D\xi \\
    \nonumber
    &=& \left(1-\tfrac12\tanh\tfrac{x}2\right )^{a} \left(1+\tfrac12\tanh\tfrac{x}2\right)^{b} p_n^{(\alpha,\beta)}\left(\tfrac12\tanh\tfrac{x}2\right)
    \end{eqnarray}
where $\alpha = 2a - 1$, $\beta = 2b-1$ and $p_n^{(\alpha,\beta)}$ is the $n$th Jacobi polynomial, normalized with respect to the measure \eqref{eqn:Jacobimeasure}. It turns out that this relationship between continuous Hahn polynomials and Jacobi polynomials generalizes a little-known identity due to Ramanujan \cite{ramanujan15sdi}, namely
\begin{displaymath}
  \int_{-\infty}^\infty |\Gamma(a+\I \xi)|^2 \E^{\I x\xi}\D \xi=\frac{\sqrt{\pi}\,\Gamma(a)\Gamma(a+\frac12)}{\cosh^{2a}\left(\frac{x}{2}\right)},\qquad a>0.
\end{displaymath}

Instead of `transformed continuous Hahn functions', we call these functions the \emph{tanh--Jacobi} functions \cite{iserles20fast}. It is convenient to map $\tfrac{x}{2}\rightarrow x$, both for aesthetic reasons and because this facilitates the computation of expansion coefficients in line with Subsection~4.1. Thus, in place of \eqref{contHahn_varphi}, we have
\begin{equation}
  \label{tanhJacobi}
  \varphi_n^{a,b}(x)=(1-\tanh x)^a(1+\tanh x)^b p_n^{(2a-1,2b-1)}(\tanh x),\qquad n\in\mathbb{Z}_+.
\end{equation}

\begin{figure}[tb]
  \begin{center}
    \includegraphics[width=220pt]{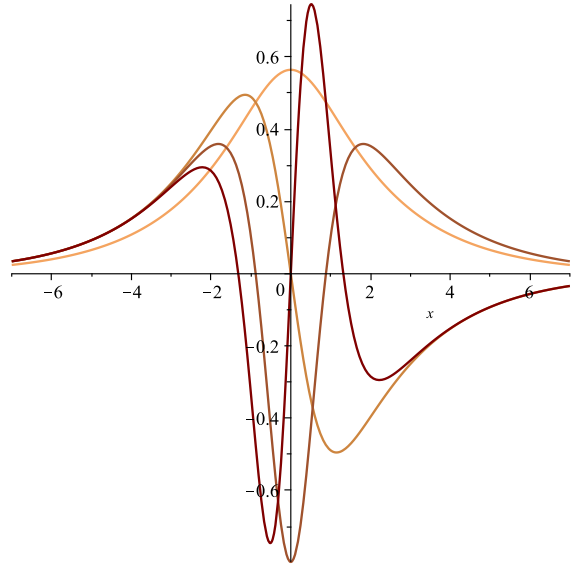}
    \caption{The tanh-Chebyshev functions $\varphi_n^{3/4,3/4}$ for $n=0,1,2,3$: darker shade corresponds to higher $n$.}
    \label{Fig:tanhChebyshev}
  \end{center}
\end{figure}

\section{Computational considerations}\label{sec:computational}
  
\subsection{Computation of expansion coefficients}
  
A major consideration in the choice of a practical basis $\Phi$ in the context of spectral methods for PDEs is the speed and ease of the calculation of the first $N$ expansion coefficients $\hat{f}_n$ such that $f(x)=\sum_n \hat{f}_n \varphi_n(x)$ for $f\in\mathrm{L}_2(\mathbb{R})$. The simplest algorithm, suitable for all bases $\Phi$, is to compute the coefficients with an $N$-point Gauss--Hermite quadrature. Provided that quadrature nodes and weights are tabulated in advance to allow for repeated calculation, each such computation requires ${\mathcal O}(N^2)$ operations if done naively. However, in the case of Hermite functions one can exploit the three-term recurrence relation for Hermite polynomials to derive an ${\mathcal O}(N\log^2N)$ algorithm \cite{leibon08fht}. The issue, though, is the stability of this procedure. Three-term recurrence relations tend to be unstable in their numerical implementation. While they can be computed stably in a bounded interval using the Clenshaw algorithm \cite{fox68cpn}, once  polynomials are orthogonal in $\mathbb{R}$, it is an elementary consequence of standard theory of orthogonal polynomials that their recurrence relations are unbounded (cf.\ for example \eqref{HF_phi}) -- this presents rather delicate implementation issues. 

At least five kinds of systems $\Phi$ in the present framework can computed faster (and stably). Letting $\E^{\I\theta}=(1+2\I x)/(1-2\I x)$ in the expression for the coefficients for Malmquist--Takenaka functions \eqref{MalmTak}, we obtain
\begin{displaymath}
  \hat{f}_n=\int_{-\infty}^\infty f(x) \overline{\varphi_n(x)}\D x=\frac{(-\I)^n}{\sqrt{2\pi}} \int_{-\pi}^\pi \left(1-\I\tan\frac{\theta}{2}\right) f\!\left(\frac12\tan\frac{\theta}{2}\right) \E^{-\I n\theta}\D\theta,
\end{displaymath}
which can be computed for $-N/2+1\leq n\leq N/2$ by FFT in $\mathcal{O}(N\log N)$ operations. The set of {\em all\/} bases $\Phi$ with this feature (namely that the expansion coefficients are equal by a monotone change of variables to the Fourier expansion coefficients of a modified function) is mildly larger: 
\begin{displaymath}
  \Phi=\left\{\gamma_n\sqrt{\frac{|\mathrm{Im}\,\lambda|}{\pi}} \E^{\I\omega x}\frac{(\lambda-x)^{n+\delta}}{(\bar{\lambda}-x)^{n+\delta+1}}\,:\, n\in\mathbb{Z}\right\},
\end{displaymath}
where $\delta,\omega\in\mathbb{R}$, $\lambda\in\mathbb{C}\setminus\mathbb{R}$ and $\gamma_n\in\mathbb{C}$, $|\gamma_n|=1$, for all $n\in\mathbb{Z}$ \cite{iserles20for}. Malmquist--Takenaka corresponds to $\delta=\omega=0$, $\lambda=\I/2$ and $\gamma_n=(-\I)^n$. 

Four systems \eqref{contHahn_varphi}, corresponding to continuous Hahn polynomials, can be computed in $\mathcal{O}(N\log N)$ operations using Fast Cosine Transform -- they correspond to 
\begin{displaymath}
  (a,b)\in\left\{(\tfrac14,\tfrac14),(\tfrac14,\tfrac34),(\tfrac34,\tfrac14),(\tfrac34,\tfrac34)\right\},
\end{displaymath}
whereby the Jacobi polynomials become Chebyshev polynomials of one of four kinds. This follows from \eqref{tanhJacobi} by the change of variables $y=\tanh x$ in the integral expression for the coefficients. 

We note in passing another approach toward the calculation of expansion coefficients. The expansion coefficients of $f$ with respect to a basis $\Phi$ coincide with the expansion coefficients of its Fourier transform with respect to the basis $P$ of orthogonal polynomials and the inner product induced by their measure. Since computing a Fourier transform with FFT costs $\mathcal{O}(N\log N)$ operations for the first $N$ coefficients, we can complement it by a  computation of conventional orthogonal expansion. While such $\mathcal{O}(N\log N)$ expansions do exist \cite{olver2020fau}, they are unfortunately restricted to Jacobi polynomials. Unless we wish to expand $f$ in PW${}_{(-1,1)}(\mathbb{R})$, this approach -- at any rate, in our current state of knowledge -- is not competitive.

\subsection{Approximation theory on the real line}

Approximation theory of analytic functions by orthogonal bases is well established in compact intervals, not so in $\mathbb{R}$. The analyticity of $f$ is assumed in a {\em Bernstein ellipse\/} surrounding an interval, whereby exponentially-fast convergence of partial sums in the underlying $\mathrm{L}_2$ norm is established, at a speed dependent on the eccentricity of ellipse.  This construction fails once the interval is infinite and no alternative overarching theory is available.

The one alternative is the classical method of steepest descent which, with a significant extent of algebraic manipulation, allows the computation of the rate of decay of expansion coefficients. The snag, though, is that each new function calls for new analysis and no general theory is available. And the little we know is baffling!

Take the Malmquist--Takenaka basis, for example. Weideman computed the rate of decay of the coefficients $\hat{f}_n$, $n\in\mathbb{Z}$, for several choices of an analytic $f$ \cite{weideman95tao}, finding
\begin{displaymath}
  \mbox{for}\;\; f(x)=\frac{1}{1+x^4}\;\;\mbox{we have}\;\; \hat{f}_n=\mathcal{O}(\rho^{-|n|}),\;\rho=1+\sqrt{2}
\end{displaymath}
--- an exponential rate of decay. Seems like the spectral decay cherished by numerical analysts. Yet,
\begin{displaymath}
  \mbox{for}\;\; f(x)=\frac{\sin x}{1+x^4}\;\;\mbox{we have}\;\; \hat{f}_n=\mathcal{O}(|n|^{-9/4})
\end{displaymath}
and the convergence is, horrifyingly, little better than quadratic. Thinking naively, $\sin x$ is an entire function, uniformly bounded in magnitude in $\mathbb{R}$: what can go wrong? The cause of the collapse in the speed of convergence is that $\sin x$ has an essential singularity at $\infty$, the North Pole of the Riemann sphere. 

Yet, the rules underlying an essential singularity at $\infty$ are hazy as well. If instead of $\sin x/(1+x^4)$ we consider $\sin x/(1+x^2)$, the rate of decay drops to $\mathcal{O}(|n|^{-5/4})$, while
\begin{eqnarray*}
  &&\mbox{for}\;\;f(x)=\E^{-x^2}\;\;\mbox{we have}\;\;\hat{f}_n=\mathcal{O}(\E^{-3|n|^{2/3}/2})\\
  \mbox{and}&&\mbox{for}\;\;f(x)=\frac{1}{\cosh x}\;\;\mbox{we have}\;\;\hat{f}_n=\mathcal{O}(\E^{-2|n|^{1/2}}).
\end{eqnarray*}
A comprehensive convergence theory for analytic functions on the real line is a significant challenge for approximation theory.

A specific type of functions of significant interest in computational quantum mechanics are {\em wave packets,\/} because in the Born--Oppenheimer formulation a wave function of a quantum system can be approximated to high accuracy by a linear combination of such functions. Thus, once we contemplate using a basis $\Phi$ as the `engine' of a spectral method to discretize PDEs of quantum mechanics, a natural question is how well it does in approximating wave packets. In a univariate setting a wave packet has the form $\cos(\omega x) \E^{-\alpha(x-x_0)^2}$, where $\alpha>0$, $\omega,x_0\in\mathbb{R}$ and typically $|\omega|\gg1$. A forthcoming paper \cite{iserles2021awp} analyses, using the method of steepest descent, the performance of different $\Phi$s in this context. It turns out that, once we wish to attain given accuracy, the performance is different for distinct orthonormal systems and that Malmquist--Takenaka functions appear to display the fastest convergence in the large $\omega$ regime.
  
  \section{Periodic bases arising from discrete orthogonal polynomials}
  Let $\mathcal{Z}\subseteq\Z$ be a set of infinite cardinality. Define the discrete inner product,
  \begin{equation}
    \label{periodic1}
    \langle f,h\rangle=\sum_{k\in\mathcal{Z}} \sigma_k f(k) \overline{h(k)},
  \end{equation}
  where $\sigma_k>0$, $k\in\mathcal{Z}$, normalized so that $\langle 1,1\rangle=1$. Typically we'll choose $\mathcal{Z}=\Z_+$ or $\Z$. The expression in \eqref{periodic1} defines an inner product, hence we can form a corresponding orthonormal polynomial system, $P=\{p_n\}_{n\in\Z_+}$.  
  
  Consider the \emph{$2\pi$-periodic} functions $\Phi = \{\varphi_n\}_{n\in\Z_+}$ given by,
  \begin{displaymath}
    \varphi_n(x)= \I^n \sum_{k\in\mathcal{Z}} \sqrt{\sigma_k} p_n(k) \E^{\I kx},\qquad n\in\Z_+,\quad x\in\R.
  \end{displaymath}
Then, first of all, $P$ being orthonormal with respect to the positive measure $\sum_{k}\sigma_k \delta(x-k)$, there exist real coefficients $B = \{b_n\}_{n\in\Z_+}$ and $C = \{c_n\}_{n\in\Z_+}$ such that
  \begin{displaymath}
    \xi p_n(\xi) =  b_{n-1}p_{n-1}(\xi) + c_np_n(\xi) + b_n p_{n+1}(\xi),\qquad n\in\Z_+,
  \end{displaymath}
  where $b_{-1}=0$ and $b_n>0$, $n\in\Z_+$. Differentiating this Fourier series term by term reveals that
  \begin{displaymath}
    \varphi_n'(x)=-b_{n-1}\varphi_{n-1}(x)+ c_n\I \varphi_n(x) + b_n\varphi_{n+1}(x),\qquad n\in\Z_+.
  \end{displaymath}
  
  Furthermore, Parseval's identity for a Fourier series gives us
  \begin{eqnarray*}
    \frac{1}{2\pi}\int_{-\pi}^\pi \varphi_m(x) \overline{\varphi_n(\xi)}  \,\D x = \I^{m-n}\sum_{k \in \mathcal{Z}} \sigma_k p_m(k)p_n(k) = \delta_{m,n}.
    \end{eqnarray*}
  Therefore, these functions are orthonormal on $\mathrm{L}_2(-\pi,\pi)$. As can be seen, Theorem \ref{thm:nonorthogonalPhi} and Theorem \ref{thm:orthogonalPhi} appear to generalize naturally to periodic functions in $\mathrm{L}_2(-\pi,\pi)$ via \emph{discrete} orthogonal polynomials. The analogue of Theorem 3 here also holds.
  
  It is of course legitimate to claim that we already have the perfect orthonormal system in $\mathrm{L}_2(-\pi,\pi) \cap C^{\infty}_{\mathrm{per}}(-\pi,\pi)$ -- the Fourier basis -- which has not a tridiagonal differentiation matrix, but a diagonal one. Of course, it also enjoys the added advantage of fast computation with FFT. Thus, it might well be that, in the greater scheme of things, new orthonormal systems of this kind are not of an immediate use. Having said so, at this stage we simply don't know!
  
  \begin{figure}[tb]
  \begin{center}
    \includegraphics[width=220pt]{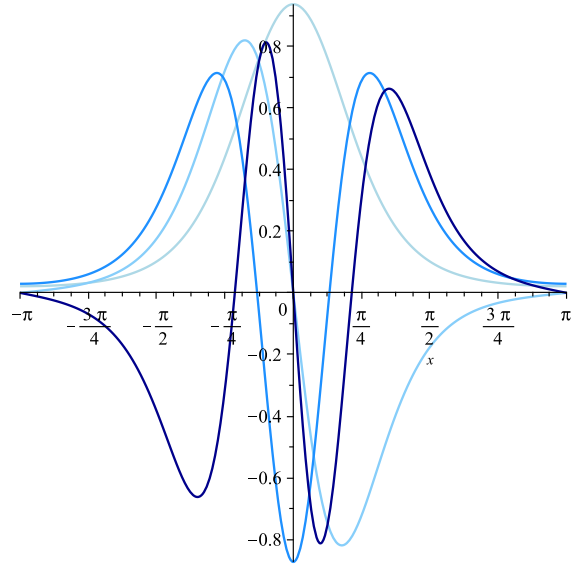}
    \caption{Transformed bilateral Charlier functions $\varphi_n$ for $n=0,1,2,3$ in $[-\pi,\pi]$: darker shade corresponds to higher $n$.}
    \label{Fig:BCharlier}
  \end{center}
\end{figure}
 
Fig.~\ref{Fig:BCharlier} displays the first four transformed functions of the bilateral Charlier measure $\D\mu_{1/2}$, where
\begin{displaymath}
  \int_{-\infty}^\infty f(\xi)\D\mu_a(\xi)=\sum_{k=-\infty}^\infty \frac{a^{|k|}}{|k|!}f(k),\qquad a>-0.
\end{displaymath}
(We symmetrize the standard Charlier measure so that the transformed functions are real.) The functions are displayed within a single period.
  
\section{Challenges and outlook}
  
This paper describes a theory in the making. The results we have uncovered so far are surprisingly elegant, perhaps because they combine the beautiful, rich theories of orthogonal polynomials and Fourier analysis. Standard orthogonal polynomials exhibit deep mathematical structure, and a major challenge is to explore if -- and how -- this structure is inherited by orthonormal systems $\Phi$.

\subsection{Transform pairs}

One thing that is particularly interesting is how in certain cases the canonical transformed functions $\Phi$ of a family of orthonormal polynomials $P$ can be expressed in terms of another family of orthonormal polynomials $Q = \{q_n\}_{n\in\Z_+}$ or in terms of known special functions. We summarize the known relationships in the following table.

\begin{table}[h]
  \begin{tabular}{l|l}
    Polynomials $p_n$  & Special functions associated to $\varphi_n$ \\[4pt]\hline
    Hermite            &  Hermite functions/polynomials               \\[2pt]
    Laguerre & Malmquist--Takenaka functions, Fourier basis, Chebyshev polynomials \\[2pt]
    Generalized Laguerre & Szeg\H{o}--Askey polynomials on the unit circle, Jacobi polynomials \\[2pt]
    Ultraspherical & Bessel functions \\[2pt]
    Continuous Hahn & Jacobi polynomials
  \end{tabular}
\end{table}

The relationship between ultraspherical polynomials and Bessel functions {\em via\/} the Fourier transform is well known to those well-versed in special functions. Less so is the relationship between generalized Laguerre polynomials and Szeg\H{o}--Askey polynomials on the unit circle as in \eqref{FL}, which appears to be noted first by the present authors in \cite{iserles20for}. The relationship between continuous Hahn polynomials and Jacobi polynomials as in \eqref{contHahn_varphi} appears not to have been noted before, but the relationship between Carlitz polynomials (which under a change of variables are a subset of continuous Hahn polynomials) and Jacobi polynomials was noted by the present authors in \cite{iserles20fast}.

\subsection{Location of zeros}

It is well known that zeros of orthogonal polynomials are real, reside in the support of the measure and that zeros of $p_{n-1}$ and $p_n$ interlace. None of these features can be taken for granted for orthonormal bases $\Phi$. Often they are -- definitely, and by design, in the case of Hermite  and continuous Hahn measures, because each $\varphi$ is a multiple of $p_n$, possibly with a strictly monotone change of argument. Sometimes the problem makes no sense: once $\D\mu$ is not symmetric with respect to the origin, $\varphi_n$ is in general complex-valued. Malmquist--Takenaka functions \eqref{MalmTak}, for example, have no real zeros at all. 

\begin{figure}[tb]
\begin{center}
  \includegraphics[width=100pt]{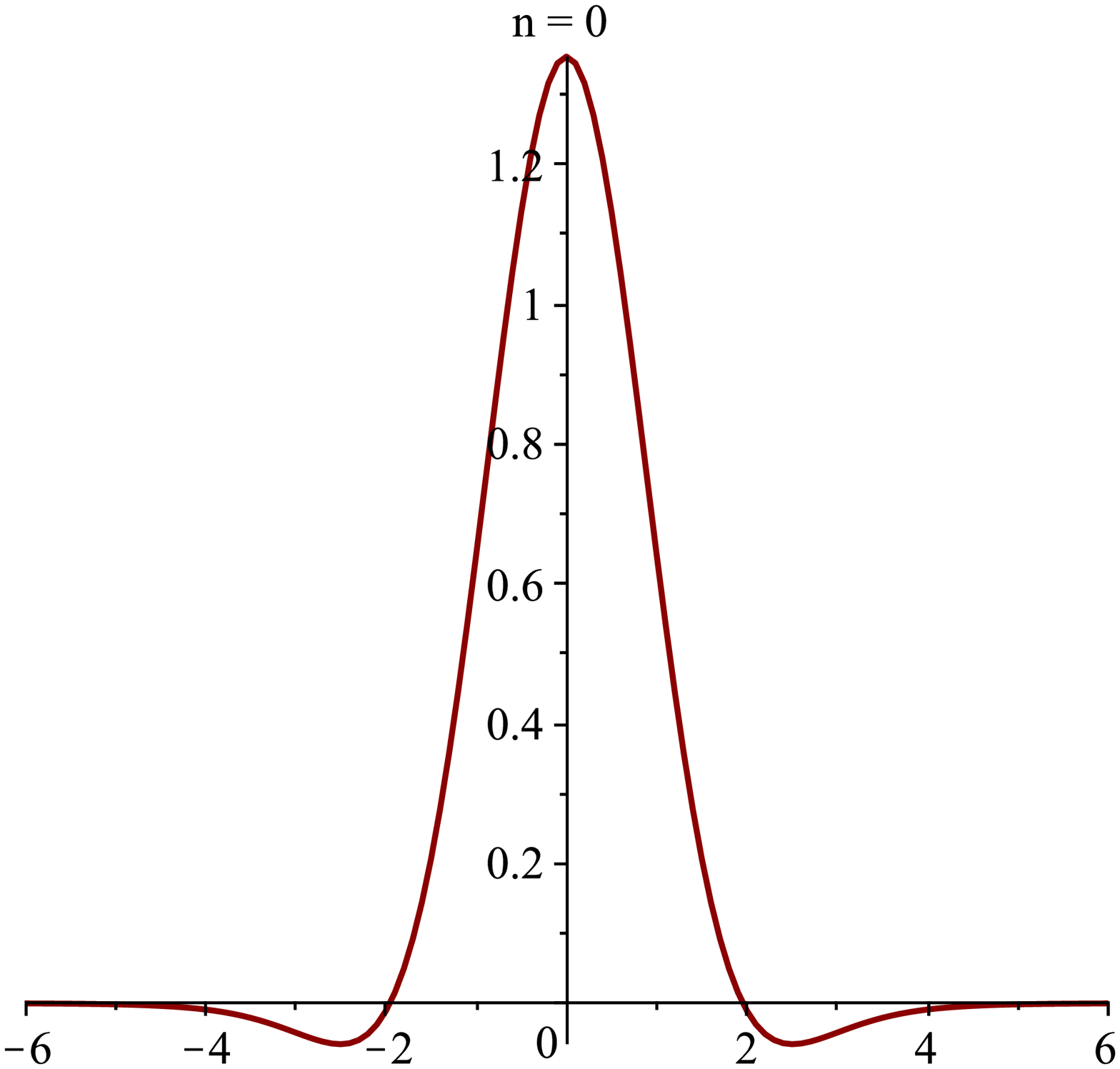}\hspace*{15pt}\includegraphics[width=100pt]{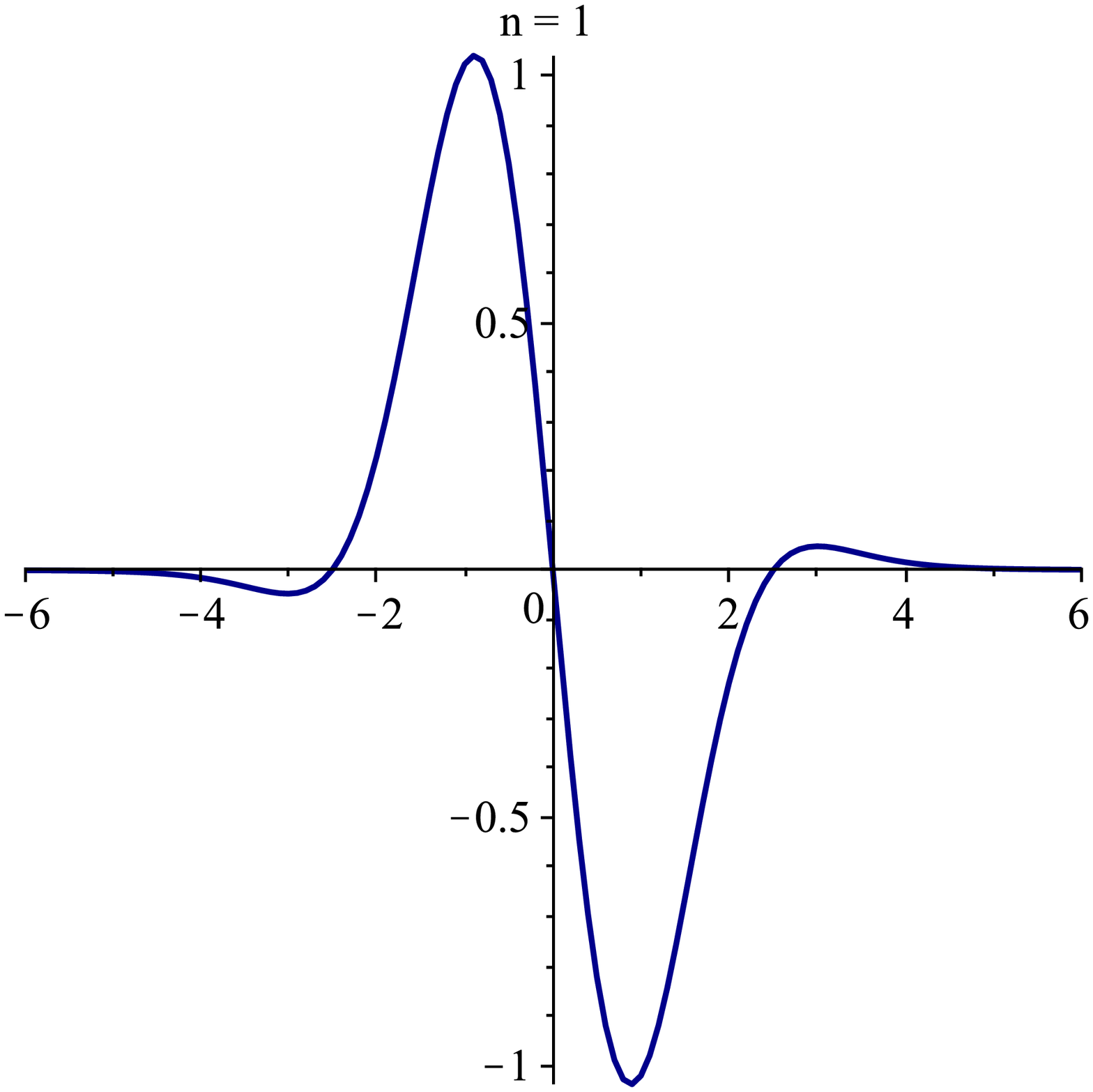}\hspace*{15pt}\includegraphics[width=100pt]{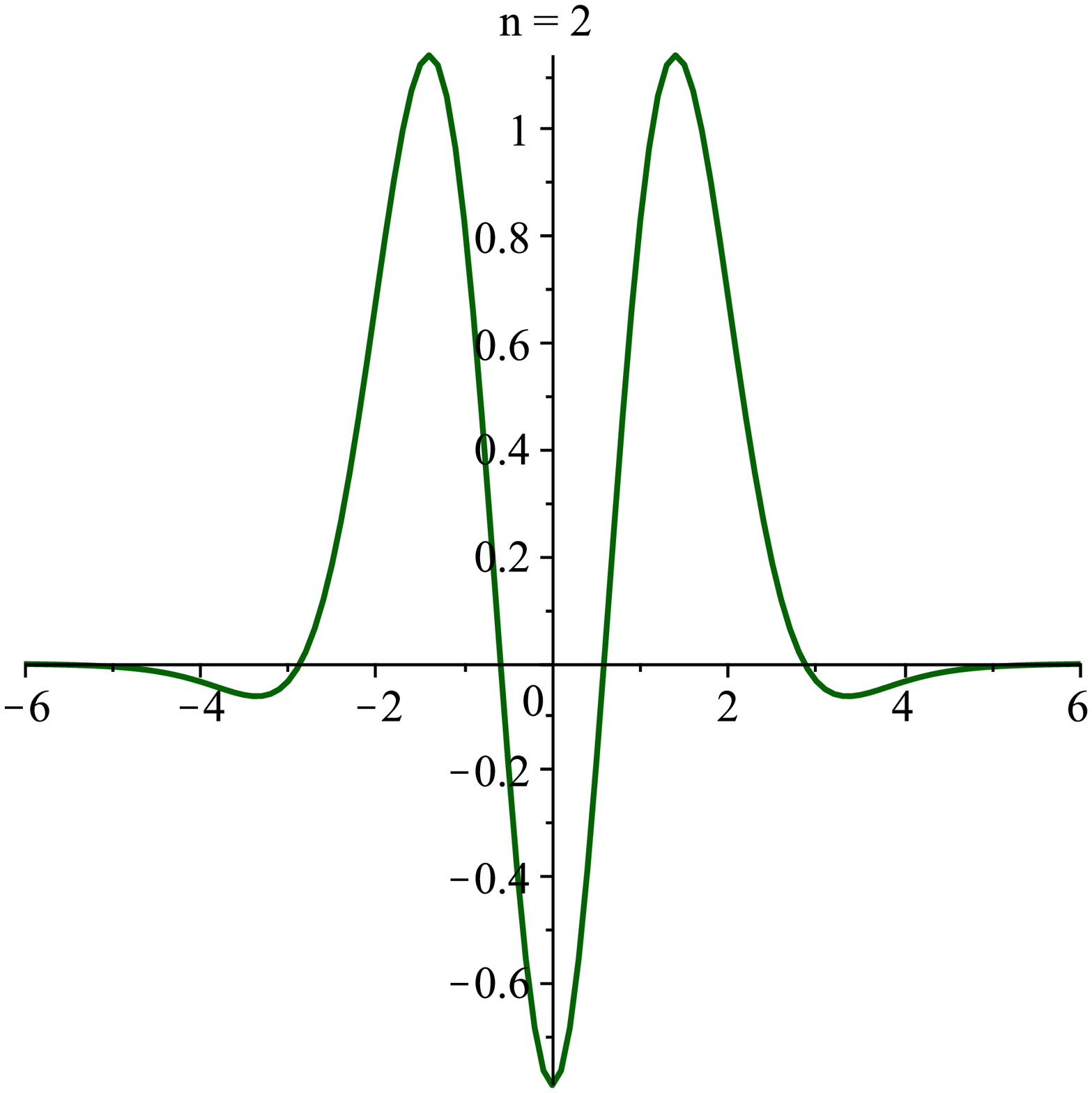}
  
  \includegraphics[width=100pt]{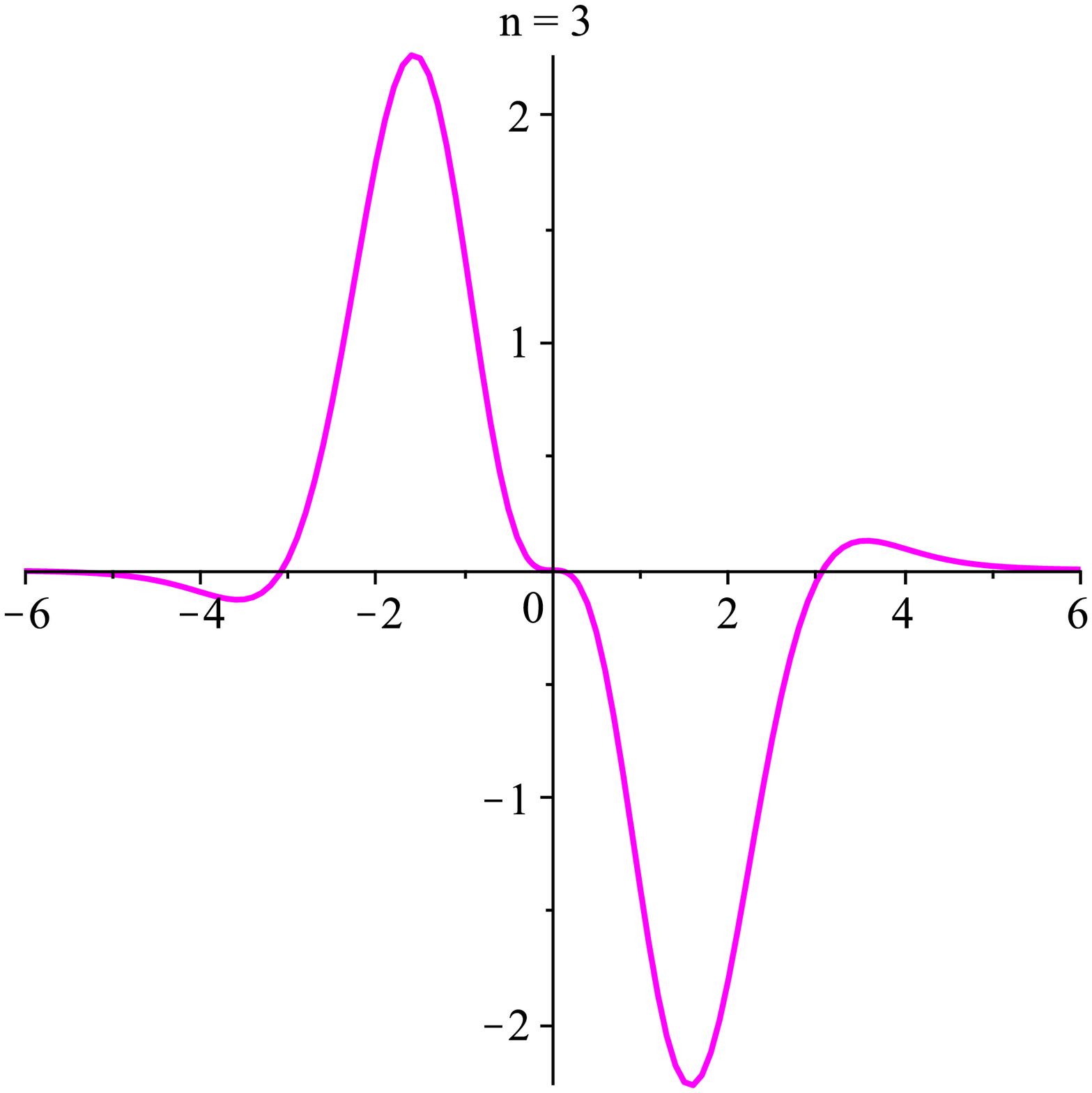}\hspace*{15pt}\includegraphics[width=100pt]{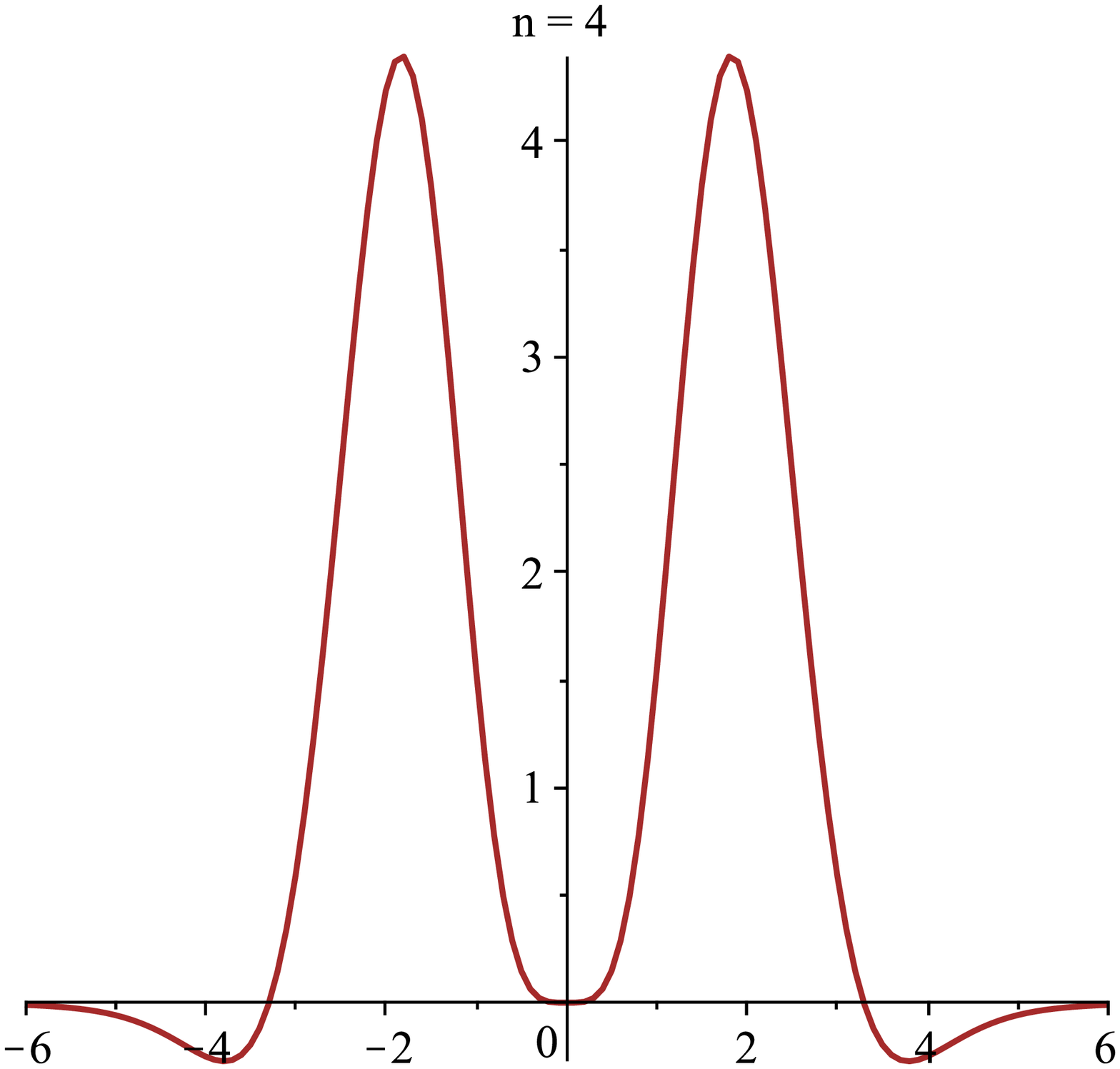}
\caption{The functions $\varphi_n$, $n=0,\ldots,4$, for $\D\mu(\xi)=(1+\xi^2)\E^{-\xi^2}\D\xi$.}
\label{Fig:SobOrtho}
\end{center}
\end{figure}

In general the picture is more hazy. Transformed Legendre functions \eqref{eqn:Legendre} have an infinity of real zeros and the zeros of $\varphi_{n-1}$ and $\varphi_n$ interlace: this follows from familiar properties of Bessel functions. The case $\D\mu(\xi)=(1+\xi^2)\E^{-\xi^2}\D\xi$ is displayed in Fig.~\ref{Fig:SobOrtho} and all  bets are off! $\varphi_0$ has two zeros, $\varphi_1$ has three and $\varphi_2$ four -- but, lest a pattern is discerned, $\varphi_3$ has just three, except that the zero at the origin has nontrivial multiplicity, while $\varphi_4$ also has three zeros while the multiplicity of the zero at the origin appears to grow. The analysis of zeros and their locations for general systems $\Phi$ is currently an open problem. 

\subsection{Sobolev orthogonality}

A natural question to ask is on the extension of our theory to more `exotic' kinds of orthogonality, e.g.\ the Sobolev $\mathrm{H}^1(\mathbb{R})$ inner product
\begin{displaymath}
  \langle f_1,f_2\rangle=\int_{-\infty}^\infty f_1(x)f_2(x)\D x+\int_{-\infty}^\infty f_1'(x)f_2'(x)\D x.
\end{displaymath}
This will be a subject of a forthcoming paper by the current authors. Spectral methods which are built of such a basis are stable in the $\mathrm{H}^1_2(\R)$ Sobolev norm, as opposed to the $\mathrm{L}_2(\R)$ norm.

\subsection{Beyond the canonical form}

The canonical form given in equation \eqref{eqn:canonical} guarantees that $\Phi$ is an orthonormal set in $\mathrm{L}_2(\R)$ and $b_n > 0$. However, this choice is not unique. Consider a basis of the form
\begin{equation}
  \label{sigma}
  \varphi_n(x) = \frac{\I^n}{\sqrt{2\pi}} \int_{-\infty}^\infty \E^{\I x \xi} p_n(\xi) \E^{\I \sigma(\xi)} |g(\xi)|^{\tfrac12} \,\D\xi,
  \end{equation}
for a measurable function $\sigma : \mathrm{supp}(\mu) \to \R$. These bases are all orthonormal in $\mathrm{L}_2(\R)$ by Theorem \ref{thm:orthogonalPhi}, and it is readily checked that all of these bases have precisely the same differentiation matrix i.e.\ the coefficients $b_n$ and $c_n$ do not depend on the function  $\sigma$. This is a subtle -- yet crucial -- distinction between orthogonal polynomials and their transformed functions. Orthogonal monic polynomials are defined uniquely by the Jacobi matrix, while a differentiation matrix is insufficient to define $\Phi$: we also need  to specify $\varphi_0$, say, yet not every $\varphi_0$ corresponds to an orthonormal system!  The choice of $\sigma$ in \eqref{sigma} captures this added freedom, while ensuring that $\Phi$ is orthonormal.

This fact has particularly interesting consequences when placed in the context of solving unitary PDEs such as Schr\"odinger's equation. A Schr\"odinger equation in one space dimension reads
\begin{displaymath}
  \frac{\partial u}{\partial t}= \I \frac{\partial^2u}{\partial x^2}+\I F(x,u),
\end{displaymath}
where $F$ is the interaction potential. There are good phenomenological reasons to solve it along the entire real line, with the initial condition $u(x,0)=u_0(x)$, $x\in\mathbb{R}$, where $u_0\in\mathrm{L}_2(\mathbb{R})$. A powerful approach toward the numerical solution of this equation is the concept of {\em splittings:\/} the solution is represented as a composition of solutions of the equation
\begin{equation}
  \label{freeS}
  \frac{\partial u}{\partial t}=\I \frac{\partial^2u}{\partial x^2}
\end{equation}
(the {\em free Schr\"odinger equation\/}) and of the ordinary differential equation $\partial u/\partial t=\I F(x,u)$ \cite{faou12gni}. This results in powerful numerical methods that recover many qualitative attributes of the solution. Suppose that we are using a spectral method with the basis $\Phi$, consistent with the theory of Section~2 and with a skew-symmetric (or skew-Hermitian) differentiation matrix. Then
\begin{displaymath}
  u_0(x)=\sum_{n\in\mathcal{Z}} \hat{u}_n \varphi_n(x)\qquad\Rightarrow\qquad u(x,t)=\sum_{n\in\mathcal{Z}} \hat{u}_n \psi_n(x,t),
\end{displaymath}
where
\begin{equation*}
  \psi_n(x,t) = \frac{\I^n}{\sqrt{2\pi}} \int_{-\infty}^\infty \E^{\I x \xi} p_n(\xi) \E^{\I \xi t^2} |g(\xi)|^{\tfrac12} \,\D\xi.
  \end{equation*}
Each $\psi_n$ is the solution of \eqref{freeS} with the initial condition $u(x,0)=\varphi_n(x)$, $x\in\mathbb{R}$ \cite{iserles2021sls}.

\subsection{A Freudian slip -- why we need more polynomials}

The size of the differentiation matrix matters! Numerous applications of spectral methods require the evaluation of the matrix exponential $\E^{\tau D_N}$ or $\E^{\I\tau D_N^2}$, where $D_N$ is the $N\times N$ principal section of $D$ and $\tau$ is the time step. Once $D$ is skew-Hermitian, these exponentials are unitary: this helps to guarantee stability in the sense of Lax.  Practical considerations, in particular the wish to use large $\tau$, require the coefficients of $D$ to be as small as possible (more specifically, the spectral radius of $D_N$ should be minimized). However, Section~2 and the connection between the entries of $D$ and recurrence coefficients of $P$, imply that the off-diagonal entries cannot be bounded: $\limsup_{n\rightarrow\infty}\alpha_n,\limsup_{n\rightarrow\infty}\gamma_n=+\infty$. Thus, for Hermite $\alpha_n=\sqrt{(n+1)/2}$ and for Malmquist--Takenaka $\alpha_n=n+1$.

The objective, thus, is to identify $P$ supported in $\mathbb{R}$ and, for simplicity, with symmetric $\D\mu$, such that the $\alpha_n$s increase at a lower rate. (Symmetry implies $\beta_n\equiv0$ and $\gamma_n=-\alpha_n$ in \eqref{diffmat}.) The obvious recourse is to use {\em Freud polynomials\/}, orthogonal with respect to $\E^{-|x|^\sigma}$, $\sigma>0$, $x\in\mathbb{R}$. According to the celebrated {\em Freud conjecture,\/}, proved by Lubinsky, Mhaskar and Saff \cite{lubinsky88pfc}, it is true that $\alpha_n=\mathcal{O}(n^{1/\sigma})$. The larger $\sigma$, the slower the decay! Yet, it is not enough to specify a measure, we also need to have the polynomials in an explicit form in order to construct the set $\Phi$, or at the very least know their recurrence coefficients. Unfortunately, and with the exception of Hermite polynomials ($\sigma=2$), the explicit form of Freud polynomials is unknown! Although there are known string relations for their recurrence coefficients, they are highly unstable as a numerical means to derive recurrence coefficients.  Thus, the challenge is to find orthogonal polynomials supported on the real line (and ideally with a symmetric measure) whose recurrence relations are known explicitly and increase slowly. 

The current state of our knowledge of orthogonal polynomials in $\mathrm{L}_2(\mathbb{R})$ is highly incomplete. Few families are known explicitly and their recurrence coefficients grow too fast for our liking. We need more polynomials!

  \bibliographystyle{plain}
  \bibliography{main}

\begin{thebibliography}{10}

\bibitem{aheizer1965classical}
Naum~Ilji{\v{c}} Aheizer and N~Kemmer.
\newblock {\em The classical moment problem and some related questions in
  analysis}.
\newblock Oliver \& Boyd Edinburgh, 1965.

\bibitem{chihara78iop}
T.~S. Chihara.
\newblock {\em An {I}ntroduction to {O}rthogonal {P}olynomials}.
\newblock Gordon and Breach Science Publishers, New York--London--Paris, 1978.
\newblock Mathematics and its Applications, Vol. 13.

\bibitem{faou12gni}
Erwan Faou.
\newblock {\em Geometric numerical integration and {S}chr\"{o}dinger
  equations}.
\newblock Zurich Lectures in Advanced Mathematics. European Mathematical
  Society (EMS), Z\"{u}rich, 2012.

\bibitem{favard35sur}
J.~Favard.
\newblock Sur les polynomes de {T}chebicheff.
\newblock {\em C.R. Acad. Sci. Paris}, 200:2052--2053, 1935.

\bibitem{fox68cpn}
L.~Fox and I.~B. Parker.
\newblock {\em Chebyshev polynomials in numerical analysis}.
\newblock Oxford University Press, London-New York-Toronto, Ont., 1968.

\bibitem{iserles2021sls}
Arieh Iserles, Karolina Kropielnicka, Katharina Schratz, and Marcus Webb.
\newblock Solving the linear {S}chr\"odinger equation on the real line.
\newblock Technical Report to appear, DAMTP, University of Cambridge, 2021.

\bibitem{iserles2021awp}
Arieh Iserles, Karen Luong, and Marcus Webb.
\newblock Approximation of wave packets on the real line.
\newblock Technical Report to appear, DAMTP, University of Cambridge, 2021.

\bibitem{iserles19oss}
Arieh Iserles and Marcus Webb.
\newblock Orthogonal systems with a skew-symmetric differentiation matrix.
\newblock {\em Found. Comput. Math.}, 19(6):1191--1221, 2019.

\bibitem{iserles20for}
Arieh Iserles and Marcus Webb.
\newblock A family of orthogonal rational functions and other orthogonal
  systems with a skew-{H}ermitian differentiation matrix.
\newblock {\em J. Fourier Anal. Appl.}, 26(1):Paper No. 19, 2020.

\bibitem{iserles20fast}
Arieh Iserles and Marcus Webb.
\newblock Fast computation of orthogonal systems with a skew-symmetric
  differentiation matrix.
\newblock {\em \mbox{To appear in} Comm. Pure Appld Maths}, 2020.

\bibitem{koekoek2010hypergeometric}
Roelof Koekoek, Peter~A Lesky, and Ren{\'e}~F Swarttouw.
\newblock {\em Hypergeometric orthogonal polynomials and their q-analogues}.
\newblock Springer Science \& Business Media, 2010.

\bibitem{lasser2020cqd}
Caroline Lasser and Christian Lubich.
\newblock Computing quantum dynamics in the semiclassical regime.
\newblock {\em Acta Numerica}, 29:229--401, 2020.

\bibitem{leibon08fht}
Gregory Leibon, Daniel~N. Rockmore, Wooram Park, Robert Taintor, and Gregory~S.
  Chirikjian.
\newblock A fast {H}ermite transform.
\newblock {\em Theoret. Comput. Sci.}, 409(2):211--228, 2008.

\bibitem{lubinsky88pfc}
D.~S. Lubinsky, H.~N. Mhaskar, and E.~B. Saff.
\newblock A proof of {F}reud's conjecture for exponential weights.
\newblock {\em Constr. Approx.}, 4(1):65--83, 1988.

\bibitem{malmquist26stc}
F.~Malmquist.
\newblock Sur la d\'etermination d'une classe de fonctions analytiques par
  leurs valeurs dans un ensemble donn\'e de points.
\newblock In {\em C.R. 6i\'eme Cong. Math. Scand. (Kopenhagen, 1925)}, pages
  253--259, Copenhagen, 1926. Gjellerups.

\bibitem{DLMF}
Frank W.~J. Olver, Daniel~W. Lozier, Ronald~F. Boisvert, and Charles~W. Clark,
  editors.
\newblock {\em N{IST} {H}andbook of {M}athematical {F}unctions}.
\newblock U.S. Department of Commerce, National Institute of Standards and
  Technology, Washington, DC; Cambridge University Press, Cambridge, 2010.
\newblock With 1 CD-ROM (Windows, Macintosh and UNIX).

\bibitem{olver2020fau}
Sheehan Olver, Richard~Mika\"el Slevinsky, and Alex Townsend.
\newblock Fast algorithms using orthogonal polynomials.
\newblock {\em Acta Numerica}, 29:573--699, 2020.

\bibitem{rainville60sf}
Earl~D. Rainville.
\newblock {\em Special functions}.
\newblock The Macmillan Co., New York, 1960.

\bibitem{ramanujan15sdi}
S.~Ramanujan.
\newblock Some definite integrals connected with {G}auss's sums [{M}essenger
  {M}ath. {\bf 44} (1915), 75--85].
\newblock In {\em Collected papers of {S}rinivasa {R}amanujan}, pages 59--67.
  AMS Chelsea Publ., Providence, RI, 2000.

\bibitem{szego1939orthogonal}
Gabor Szeg\H{o}.
\newblock {\em Orthogonal polynomials}, volume~23.
\newblock American Mathematical Soc., 1939.

\bibitem{takenaka25oof}
S.~Takenaka.
\newblock On the orthogonal functions and a new formula of interpolation.
\newblock {\em Japanese J. Maths}, 2:129--145, 1926.

\bibitem{weideman95tao}
J.A.C. Weideman.
\newblock Theory and applications of an orthogonal rational basis set.
\newblock In {\em Proceedings South African Num. Math. Symp 1994, Univ. Natal},
  1994.

\end{thebibliography}

\end{document}